\documentclass[12pt,a4paper,twoside]{article}

\newcommand{\pageformat}[6]{\setlength{\hoffset}{-1in}
                  \setlength{\voffset}{-1in}
                  \addtolength{\hoffset}{#5}
                            \addtolength{\voffset}{#6}
                            \setlength{\oddsidemargin}{#1}
                            \setlength{\evensidemargin}{#2}
                            \setlength{\textwidth}{\paperwidth}
                  \addtolength{\textwidth}{-\oddsidemargin}
                  \addtolength{\textwidth}{-\evensidemargin}
                  \addtolength{\textwidth}{-\marginparsep}
                  \addtolength{\textwidth}{-\marginparwidth}
                            \setlength{\topmargin}{#3}
                            \setlength{\textheight}{\paperheight}
                  \addtolength{\textheight}{-\topmargin}
                  \addtolength{\textheight}{-\headheight}
                  \addtolength{\textheight}{-\headsep}
                  \addtolength{\textheight}{-\footskip}
                  \addtolength{\textheight}{-#4}}
\pageformat{2cm}{3cm}{24mm}{24mm}{1pt}{0pt}

\usepackage{ifthen}
\newboolean{article}
    \setboolean{article}{true}
\newboolean{report}
\newboolean{book}
\newboolean{letter}
\newboolean{german}
\newboolean{italian}
\newboolean{nobaselinestretch}
\newboolean{nosectionappendix}
\newboolean{oldtoc}
\newboolean{nosectionequn}
\newboolean{notheorem}
\ifthenelse{\boolean{german}}{
    \usepackage{german}}{}

\usepackage[latin1]{inputenc}
\usepackage{amsmath}
\usepackage{amssymb}
\usepackage[mathscr]{eucal}

\ifthenelse{\boolean{notheorem}}{}{
    \usepackage{theorem}}

\ifthenelse{\boolean{nobaselinestretch}}{}{
    \renewcommand{\baselinestretch}{1.25}}

\newenvironment{env}[2]{\begin{#1}#2\end{#1}}{}
    \newcommand{\beq}[1]{\begin{env}{equation}{#1}}
    \newcommand{\beqn}[1]{\begin{env}{equation*}{#1}}
    \newcommand{\bal}[1]{\begin{env}{align}{#1}}
    \newcommand{\baln}[1]{\begin{env}{align*}{#1}}
    \newcommand{\bga}[1]{\begin{env}{gather}{#1}}
    \newcommand{\bgan}[1]{\begin{env}{gather*}{#1}}
    \newcommand{\bflal}[1]{\begin{env}{flalign}{#1}}
    \newcommand{\bflaln}[1]{\begin{env}{flalign*}{#1}}
    \newcommand{\bmu}[1]{\begin{env}{multline}{#1}}
    \newcommand{\bmun}[1]{\begin{env}{multline*}{#1}}
    \newcommand{\bsp}[1]{\begin{env}{split}{#1}}

    \newcommand{\eeq}{\end{env}}
    \newcommand{\eeqn}{\end{env}}
    \newcommand{\eal}{\end{env}}
    \newcommand{\ealn}{\end{env}}
    \newcommand{\ega}{\end{env}}
    \newcommand{\egan}{\end{env}}
    \newcommand{\eflal}{\end{env}}
    \newcommand{\eflaln}{\end{env}}
    \newcommand{\emu}{\end{env}}
    \newcommand{\emun}{\end{env}}
    \newcommand{\esp}{\end{env}}

\newcommand{\lf}{\vspace{2ex}}

\renewcommand{\bf}[1]{\textbf{#1}}
\renewcommand{\it}[1]{\textit{#1}}

\renewcommand{\sf}[1]{\textsf{#1}}

\renewcommand{\tt}[1]{\texttt{#1}}
\newcommand{\hl}[1]{\bf{\it{#1}}}

\newcommand{\mbf}[1]{\mathbf{#1}}
\newcommand{\msf}[1]{\text{\small$\sf{#1}$}}

\newcommand{\cmc}[1]{\mathcal{#1}}
\newcommand{\eus}[1]{\mathscr{#1}}
\newcommand{\euf}[1]{\mathfrak{#1}}
\newcommand{\bb}[1]{\mathbb{#1}}

\newcommand{\mfootnotesize}[1]{{\setlength{\arraycolsep}{.5ex}\text{\footnotesize$#1$}}}
\newcommand{\mscriptsize}[1]{{\setlength{\arraycolsep}{.3ex}\text{\scriptsize$#1$}}}
\newcommand{\mtiny}[1]{{\setlength{\arraycolsep}{.3ex}\text{\tiny$#1$}}}
\newcommand{\nbd}[1]{$#1$\nobreakdash--}
\newcommand{\ol}[1]{\overline{#1}}
\newcommand{\ul}[1]{\underline{#1}}

\newcommand{\wh}[1]{\widehat{#1}}

\newcommand{\vt}{\vartheta}

\newcommand{\norm}[1]{\left\lVert#1\right\rVert}

\newcommand{\family}[1]{\left(#1\right)}

\newcommand{\bfam}[1]{\bigl(#1\bigr)}
\newcommand{\Bfam}[1]{\Bigl(#1\Bigr)}
\newcommand{\AB}[1]{\langle#1\rangle}
\newcommand{\bAB}[1]{\bigl\langle#1\bigr\rangle}

\newcommand{\CB}[1]{\{#1\}}
\newcommand{\bCB}[1]{\bigl\{#1\bigr\}}
\newcommand{\BCB}[1]{\Bigl\{#1\Bigr\}}
\newcommand{\SB}[1]{[#1]}

\newcommand{\RO}[1]{[#1)}

\newcommand{\Matrix}[1]{\begin{pmatrix}#1\end{pmatrix}}

\newcommand{\fMatrix}[1]{\mfootnotesize{\Matrix{#1}}}
\newcommand{\sMatrix}[1]{\mscriptsize{\Matrix{#1}}}
\newcommand{\tMatrix}[1]{\mtiny{\Matrix{#1}}}

\newcommand{\set}[2][]{
    \ifthenelse{\equal{#1}{}}{
        \CB{#2}}{
        \CB{#1~|~#2}}}
\newcommand{\bset}[2][]{
    \ifthenelse{\equal{#1}{}}{
        \bCB{#2}}{
        \bCB{#1~|~#2}}}
\newcommand{\Bset}[2][]{
    \ifthenelse{\equal{#1}{}}{
        \BCB{#2}}{
        \BCB{#1~\big|~#2}}}
\newcommand{\zero}{\CB{0}}

\DeclareMathOperator{\ls}{\normalfont\msf{span}}
\DeclareMathOperator{\cls}{\ol{\ls}}

\DeclareMathOperator{\id}{\normalfont\msf{id}}

\newcommand{\C}{\bb{C}}

\newcommand{\E}{\bb{E}}

\newcommand{\N}{\bb{N}}

\newcommand{\R}{\bb{R}}
\newcommand{\bS}{\bb{S}}

\newcommand{\cA}{\cmc{A}}
\newcommand{\cB}{\cmc{B}}
\newcommand{\cC}{\cmc{C}}

\newcommand{\cK}{\cmc{K}}

\newcommand{\sB}{\eus{B}}

\newcommand{\sF}{\eus{F}}

\newcommand{\sK}{\eus{K}}

\newcommand{\sS}{\eus{S}}

\newcommand{\en}{\euf{n}}

\newcommand{\eu}{\euf{u}}

\newcommand{\eH}{\euf{H}}

\newcommand{\U}{\mbf{1}}

\newcommand{\G}{\Gamma}

\newcommand{\I}{{I\!\!\!\;I}}

\newcommand{\s}{\text{\scriptsize$\sS$}}

\ifthenelse{\boolean{nosectionequn}}{}{
    \numberwithin{equation}{section}
    }

\ifthenelse{\boolean{article}\or\boolean{letter}\or\boolean{nosectionequn}}{
    \setboolean{nosectionappendix}{true}}{}
\ifthenelse{\boolean{nosectionappendix}}{}{
    \renewcommand{\appendix}{
        \chapter*{\appendixname}
        \addcontentsline{toc}{chapter}{\appendixname}
        \renewcommand{\thesection}{\Alph{section}}
        \setcounter{section}{0}}}
   
\ifthenelse{\boolean{report}\or\boolean{book}}{
    }{}

\ifthenelse{\boolean{notheorem}}{}{
        \newcommand{\mnname}{Mathematical note.}
        \newcommand{\enname}{End of the note.}
        \newcommand{\definame}{Definition.}
        \newcommand{\propname}{Proposition.}
        \newcommand{\lemname}{Lemma.}
        \newcommand{\exname}{Example.}
        \newcommand{\exername}{Exercise.}
        \newcommand{\remname}{Remark.}
        \newcommand{\obname}{Observation.}
        \newcommand{\thmname}{Theorem.}
        \newcommand{\corname}{Corollary.}
        \newcommand{\proofname}{Proof.}
    \ifthenelse{\boolean{german}}{
        \renewcommand{\mnname}{Mathematische Notiz.}
        \renewcommand{\enname}{Ende der Notiz.}
        \renewcommand{\exname}{Beispiel.}
        \renewcommand{\exername}{Übung.}
        \renewcommand{\remname}{Bemerkung.}
        \renewcommand{\obname}{Beobachtung.}
        \renewcommand{\thmname}{Satz.}
        \renewcommand{\corname}{Korollar.}
        \renewcommand{\proofname}{Beweis.}}{}
    \ifthenelse{\boolean{italian}}{
        \renewcommand{\mnname}{Nota matematica.}
        \renewcommand{\enname}{Fina della nota.}
        \renewcommand{\definame}{Definizione.}
        \renewcommand{\propname}{Proposizione.}
        \renewcommand{\exname}{Esempio.}
        \renewcommand{\exername}{Esercizio.}
        \renewcommand{\remname}{Nota.}
        \renewcommand{\obname}{Osservazione.}
        \renewcommand{\thmname}{Teorema.}
        \renewcommand{\corname}{Corollario.}
        \renewcommand{\proofname}{Dimostrazione.}

       \renewcommand{\appendixname}{Appendice}

       }{}
    \theoremheaderfont{\normalfont\bfseries}
    \theoremstyle{change}
        \theorembodyfont{\rmfamily}
            \newtheorem{emp}{}[section]
                \newcommand{\bemp}[1][]{
                    \begin{emp}\hskip-\labelsep\bf{#1}\hskip\labelsep}
                \newcommand{\eemp}{\end{emp}}
\newtheorem{itemp}[emp]{}
                \newcommand{\bitemp}[1][]{
                    \begin{itemp}\hskip-\labelsep\bf{#1}\hskip\labelsep\normalfont\itshape}
                \newcommand{\eitemp}{\end{itemp}}
            \newtheorem{mn}[emp]{\mnname}
                \newcommand{\bnm}{\begin{mn}~\begin{quotation}\renewcommand{\baselinestretch}{1}\small\noindent\ignorespaces}
                \newcommand{\enm}{\end{quotation}\hfill\bf{\enname}\end{mn}}
            \newtheorem{ex}[emp]{\exname}
                \newcommand{\bex}{\begin{ex}}
                \newcommand{\eex}{\end{ex}}
            \newtheorem{exer}[emp]{\exername}
                \newcommand{\bexer}{\begin{exer}}
                \newcommand{\eexer}{\end{exer}}
            \newtheorem{defi}[emp]{\definame}
                \newcommand{\bdefi}{\begin{defi}}
                \newcommand{\edefi}{\end{defi}}
            \newtheorem{rem}[emp]{\remname}
                \newcommand{\brem}{\begin{rem}}
                \newcommand{\erem}{\end{rem}}
            \newtheorem{ob}[emp]{\obname}
                \newcommand{\bob}{\begin{ob}}
                \newcommand{\eob}{\end{ob}}
        \theorembodyfont{\normalfont\itshape}
            \newtheorem{thm}[emp]{\thmname}
                \newcommand{\bthm}{\begin{thm}}
                \newcommand{\ethm}{\end{thm}}
            \newtheorem{prop}[emp]{\propname}
                \newcommand{\bprop}{\begin{prop}}
                \newcommand{\eprop}{\end{prop}}
            \newtheorem{cor}[emp]{\corname}
                \newcommand{\bcor}{\begin{cor}}
                \newcommand{\ecor}{\end{cor}}
            \newtheorem{lem}[emp]{\lemname}
                \newcommand{\blem}{\begin{lem}}
                \newcommand{\elem}{\end{lem}}
\newenvironment{empn}[1]{\lf\noindent\bf{#1}\ignorespaces\hskip\labelsep}{\lf}
		\newcommand{\bempn}[1]{\begin{empn}{#1}}
		\newcommand{\eempn}{\end{empn}}
		\newcommand{\bitempn}[1]{\begin{empn}{#1}\normalfont\itshape}
		\newcommand{\eitempn}{\end{empn}}
                \newcommand{\bnmn}{\begin{empn}{\mnname}~\begin{quotation}\renewcommand{\baselinestretch}{1}\small\noindent\ignorespaces}
                \newcommand{\enmn}{\end{quotation}\hfill\bf{\enname}\end{empn}}
		\newcommand{\bexn}{\begin{empn}{\exname}}
		\newcommand{\eexn}{\end{empn}}
		\newcommand{\bexern}{\begin{empn}{\exername}}
		\newcommand{\eexern}{\end{empn}}
		\newcommand{\bdefin}{\begin{empn}{\definame}}
		\newcommand{\edefin}{\end{empn}}
		\newcommand{\bremn}{\begin{empn}{\remname}}
		\newcommand{\eremn}{\end{empn}}
		\newcommand{\bobn}{\begin{empn}{\obname}}
		\newcommand{\eobn}{\end{empn}}

\newcommand{\qedsymbol}{~\rule[-0.35mm]{2mm}{2mm}}
    \newcounter{proof}[emp]
    \newenvironment{Proof}[1]{
        \vspace{1ex}
        \renewcommand{\item}[1][\stepcounter{proof}(\roman{proof})]%
            {##1\hskip\labelsep}
        \noindent\textsc{#1\hskip\labelsep}}{
        \nolinebreak\qedsymbol}
    \newcommand{\proof}[1][\proofname]{
        \begin{Proof}{#1}\ignorespaces}
    \newcommand{\qed}{\end{Proof}}
    \newcommand{\noqed}{
        \renewcommand{\qedsymbol}{}
        \end{Proof}}}
    \ifthenelse{\boolean{italian}}{
        \renewcommand{\proofname}{Dimostrazione.}}{}

\usepackage[varg]{txfonts}

\renewcommand{\thefootnote}{[\arabic{footnote}]}

\newcommand{\cvN}{\euf{cvN}}
\newcommand{\ecC}{\euf{cC}}

\usepackage{bbm}

\begin{document}

\title{Product Systems; a Survey with Commutants in View\renewcommand{\thefootnote}{}\thanks{2000 AMS-Subject classification: 46L53; 46L55; 46L08; 60J25; 81S25; 12H20; }}
\author{}
\author{
~\\
Michael Skeide\thanks{This work is supported by research funds of University of Molise (Dipartimento S.E.G.e S.) and Italian MIUR (PRIN 2005).}\\[1ex]
{\small\itshape Dipartimento S.E.G.e S.}\\
{\small\itshape Universit\`a\ degli Studi del Molise}\\
{\small\itshape Via de Sanctis}\\
{\small\itshape 86100 Campobasso, Italy}\\
{\small{\itshape E-mail: \tt{skeide@unimol.it}}}\\
{\small{\itshape Homepage: \tt{http://www.math.tu-cottbus.de/INSTITUT/lswas/\_skeide.html}}}\\
\\
%\sc{Under Construction}\\
}
\date{January 2007}
%$\mtt{\ul{~}}$

{
%\addtolength{\parskip}{-1ex}
\renewcommand{\baselinestretch}{1}
\maketitle

%\vfill
%\newpage

%\vspace{-3ex}

%\clearpage

\begin{abstract}
\noindent
The theory of product systems both of Hilbert spaces (Arveson systems) and product systems of Hilbert modules has reached a status where it seems appropriate to rest a moment and to have a look at what is known so far and what are open problems. However, the attempt to give an approximately complete account in view pages is destined to fail already for Arveson systems since Tsirelson, Powers and Liebscher have discovered their powerful methods to construct large classes of examples. In this survey we concentrate on that part of the theory that works also for Hilbert modules. This does not only help to make a selection among the possible topics, but it also helps to shed some new light on the case of Arveson systems. Often, proofs that work for Hilbert modules also lead to simpler proofs in the case of Hilbert spaces. We put emphasis on those aspects that arise from recent results about \it{commutants} of von Neumann correspondences, which, in the case of Hilbert spaces, explain the relation between the Arveson system and the Bhat system associated with an \nbd{E_0}semigroup on $\sB(H)$.
\end{abstract}

%\tableofcontents
}

%\vfill

%{\parskip0.5ex plus 0.5ex minus 0.5ex

\section{Introduction}\label{intro}

A \it{product system of Hilbert spaces} is a family $E^\otimes=\bfam{E_t}_{t\in\R_+}$ of Hilbert spaces that factor as
\beqn{
E_{s+t}
~\cong~
E_s\otimes E_t
}\eeqn
by means of an associative bilinear multiplication $E_s\times E_t\ni(x_s,y_t)\mapsto x_sy_t\in E_{s+t}$. (Depending on the application, there are also technical conditions about continuity or measurability of sections. We speak about this later on.) The definition of such product systems is due to Arveson \cite{Arv89}. It is motivated by Arveson's construction that associates with every \it{\nbd{E_0}semigroup} (a semigroup of normal unital endomorphisms) of the $\sB(H)$ of all adjointable operators on a Hilbert space $H$ a product system ${E^A}^\otimes$. If $H$ is infinite-dimensional and separable, then the product system determines the \nbd{E_0}semigroup up to \it{cocycle conjugacy}. In a series of four articles \cite{Arv89,Arv90a,Arv89a,Arv90} Arveson showed the \it{fundamental theorem}, namely, that every product system of Hilbert spaces is the one associated with a suitable \nbd{E_0}semigroup. Thus, there is a one-to-one correspondence between product systems (up to isomorphism) and \nbd{E_0}semigroups (up to cocycle conjugacy). In the sequel, we will speak about an \hl{Arveson system} if we intend a product system of Hilbert spaces. In particular, we will speak about the Arveson system associated with an \nbd{E_0}semigroup on $\sB(H)$.

Meanwhile, product systems of Hilbert \it{bimodules} or, more fashonably, \it{correspondences} made appearance in many contexts. Bhat and Skeide \cite{BhSk00} constructed a product system of correspondences over a (unital) \nbd{C^*}algebra $\cB$ from a (unital) CP-semigroup on $\cB$. (See also the discussion of Muhly and Solel \cite{MuSo02} in Remark \ref{finrem}.) This construction overcomes constructions by Bhat \cite{Bha96} and by Arveson \cite{Arv96} who construct an Arveson system starting from a CP-semigroup on $\sB(H)$ by, first, \it{dilating} in a unique way the CP-semigroup to a \it{minimal} \nbd{E_0}semigroup and, then, constructing the Arveson system of that \nbd{E_0}semigroup. The construction of \cite{BhSk00}, instead, is direct and allows, then, to construct the minimal dilation in a transparent way. Only later, Skeide \cite{Ske02,Ske05a,Ske04p} associated in several ways with an \nbd{E_0}semigroup a product system. Now the \nbd{E_0}semigroup acts on the algebra $\sB^a(E)$ of all adjointable maps on a Hilbert \nbd{\cB}module $E$. Although historically earlier, the approach to product systems from CP-semigroups (that is, \it{irreversible} quantum dynamics) has the disadvantage that not all product systems arise in that way. While one of the latest results (Skeide \cite{Ske07p}; still in preparation) asserts that \it{cum grano salis} every product system comes from an \nbd{E_0}semigroup (that is, \it{reversible} quantum dynamics in a sense we specify later on). So the approach via \nbd{E_0}semigroups allows a more coherent discussion. In this survey we will concentrate on this connection between product systems and \nbd{E_0}semigroups, while we will have no space to discuss also the connections with CP-semigroup and their dilations; see Skeide \cite{Ske03b}. Also basic classification of product systems must be sacrificed; see Skeide \cite{Ske03b,Ske06d}.

The basic factorization property of the symmetric Fock space
\beqn{
\G(H_1\oplus H_2)
~=~
\G(H_1)\otimes\G(H_2)
}\eeqn
($H_1$ and $H_2$ some Hilbert spaces) has drawn attention since a long time. In the form
\beqn{\tag{$*$}\label{Fockfac}
\G(L^2(\SB{r,t},K))
~=~
\G(L^2(\SB{r,s},K))\otimes\G(L^2(\SB{s,t},K)),
~~~~~~~~~~~~
r\le s\le t
}\eeqn
($K$ a Hilbert space) it made appearance in the work of Araki \cite{Ara70} and Streater \cite{Str69} on current representations of Lie algebras, in the work of Parthasarathy and Schmidt \cite{PaSchm72} about Lévy processes (culminating in Schürmann's work \cite{MSchue93} on quantum Lévy processes) and in quantum stochastic calculus on the symmetric Fock space initiated by Hudson and Parthasarathy \cite{HuPa84}.

Let us put $E_t=\G(L^2(\SB{0,t},K))$. Then, from the beginning, there are two possibilities to use \eqref{Fockfac} in order to define an isomorphism $E_s\otimes E_t\cong E_{s+t}$, namely,
\baln{
E_s\otimes E_t
&
~\cong~
\s_tE_s\otimes E_t
~\cong~
E_{s+t}
&
&
\text{and}
&
E_s\otimes E_t
&
~\cong~
E_s\otimes\s_sE_t
~\cong~
E_{s+t},
}\ealn
where $\s_t\colon\G(L^2(\SB{0,s},K))\rightarrow\G(L^2(\SB{t,t+s},K))$ is the time shift. If we consider the \it{CCR-flow}, that is, the \nbd{E_0}semigroup induced on $\sB\bfam{\G(L^2(\R_+,K))}$ by the time shift, then the associated Arveson system is $E_t$ with the second choice of an isomorphism, that is, with the time shift acting on the right factor in $E_s\otimes E_t$. However, Bhat discovered a second possibility to associate an Arveson with an \nbd{E_0}semigroup. In the case of the CCR-flow one obtains the same Hilbert spaces $E_t$ but with the first choice of an isomorphism, that is, with the time shift acting on the left factor in $E_s\otimes E_t$. More generally, the \hl{Bhat system} associated with any \nbd{E_0}semigroup shows always to be \it{anti-isomorphic} to the associated Arveson system.

This ambivalence in the tensor product of Hilbert spaces, where we may switch the order of factors without changing (up to canonical isomorphism) the resulting space, is by far less innocent than it appears at the first sight. Nothing like this is true in the module case for the tensor product of correspondences over $\cB$. (It is very well possible that in one order their tensor product is $\zero$, while in the other order it is not.) In fact, we will see that the construction of a product system of correspondences over $\cB$ from an \nbd{E_0}semigroup on $\sB^a(E)$ for some Hilbert \nbd{\cB}module $E$ corresponds to the construction of the Bhat system of an \nbd{E_0}semigroup on $\sB(H)$. Also the construction of product system following the ideas of Arveson is still possible. However, it yields a product system of correspondences over the commutant $\cB'$ of $\cB$ and works nicely only for von Neumann algebras $\cB$. The relation between these two product systems is that one is the \it{commutant} of the other. The commutant of a correspondence was introduced in Skeide \cite{Ske03c}, the conribution to the proceedings of the conference in Mount Holyoke 2002.

In the space available we are not able to even scratch the basic classification results for product systems. We refere the reader to the still quite up-to-date survey Skeide \cite{Ske03b} in the proceedings of the conference in Burg 2001. The classification is based on spatial product systems and their product in Skeide \cite{Ske06d} (preprint 2001).

\lf
We fix some notations used throughout, and recall very few basics about Hilbert modules in order to make this survey digestable also for nonexperts in Hilbert modules.

Let $\cB$ be a \nbd{C^*}algebra. Recall that a \hl{pre-Hilbert \nbd{\cB}module} is a right \nbd{\cB}module $E$ with a sesquilinear inner product $\AB{\bullet,\bullet}\colon E\times E\rightarrow\cB$ satisfying $\AB{x,x}\ge0$ for all $x\in E$ (positivity), $\AB{x,yb}=\AB{x,y}b$ for all $x,y\in E;b\in\cB$ (right linearity), and $\AB{x,x}=0~\Longrightarrow~x=0$ (definiteness). If definiteness is missing, then $E$ is a \hl{semi-Hilbert \nbd{\cB}module}. (Properties like $\AB{x,y}^*=\AB{y,x}$ and $\AB{xb,y}=b^*\AB{x,y}$ are automatic.) The most basic property of the inner product in a semi-Hilbert \nbd{\cB}module is the following \hl{Cauchy-Schwartz inequality}
\beqn{
\AB{x,y}\AB{y,x}
~\le~
\norm{\AB{y,y}}\AB{x,x}.
}\eeqn
By Cauchy-Schwartz inequality it is possible to quotient out length-zero elements. By Cauchy-Schwartz inequality $\norm{x}:=\sqrt{\AB{x,x}}$ defines a norm on the pre-Hilbert module $E$. If $E$ is complete in that norm, then $E$ is a \hl{Hilbert \nbd{\cB}module}. By Cauchy-Schwartz inequality the operator norm turns the algebra of bounded adjointable operators $\sB^a(E)$ on the pre-Hilbert module $E$ into a pre--\nbd{C^*}algebra. Recall that a map $a$ on $E$ is \hl{adjointable}, if it admits an adjoint $a^*$ such that $\AB{x,ay}=\AB{a^*x,y}$ for all $x,y\in E$. Every adjointable map is closeable. Therefore, by the \it{closed graph theorem}, an adjointable map on a Hilbert module is bounded, automatically.

In order to speak about product systems we need the (internal) tensor product, and the tensor product is among bimodules or \it{correspondences}. If $\cA$ is another \nbd{C^*}algebra, then a \hl{correspondence} from $\cA$ to $\cB$ (or a \nbd{\cA}\nbd{\cB}correspondence) is a Hilbert \nbd{\cB}module with a \hl{nondegenerate} representation of $\cA$ by adjointable operators. If $\cA=\cB$, then we speak also of a correspondence over $\cB$ or of a \nbd{\cB}correspondence.\footnote{The nondegeneracy condition is crucial in all what follows. For the right action of $\cB$ on a Hilbert \nbd{\cB}module it is automatic. (Exercise: Why?) But, there are left actions that act degenerately. However, in that case we will never say \nbd{\cA}\nbd{\cB}module, but rather speak of a (possibly degenerate) representation of $\cA$.} The (\hl{internal}) \hl{tensor product} of a correspondence $E$ from $\cA$ to $\cB$ and a correspondence $F$ from $\cB$ to $\cC$ is the unique correspondence $E\odot F$ from $\cA$ to $\cC$ that is generated by elementary tensors $x\odot y$ with inner product
\beq{\label{tpdef}
\AB{x\odot y,x'\odot y'}
~=~
\bAB{y',\AB{x,x'}y'}
}\eeq
and the obvious bimodule operation. Uniqueness is, in the sense of a universal property, up to canonical isomorphism. (In two realizations, simply identify the elementary tensors. For a construction take the vector space tensor product $E\otimes F$, define a semiinner product by \eqref{tpdef} and divide by the length-zero elements.) The tensor product applies also if $E$ is just a Hilbert \nbd{\cB}module, as every Hilbert \nbd{\cB}module $E$ may be viewed as a correspondence from $\sB^a(E)$ to $\cB$. This also shows that $E\odot F$ carries a canonical nondegenerate left action of $a\in\sB^a(E)$ which we denote by $a\odot\id_F$ or, sometimes, simply by $a$, too. (Attention! The unital embedding $\sB^a(E)\rightarrow\sB^a(E)\odot\id_F\subset\sB^a(E\odot F)$ need not be faithful.) By $\sB^{a,bil}(F)$ we denote the space of those elements $a\in\sB^a(F)$ that are \hl{bilinear}, that is, which fulfill $a(by)=b(ay)$ for all $b\in\cB,y\in F$. There is an embedding $\sB^{a,bil}(F)\rightarrow\id_E\odot\sB^{a,bil}(F)\subset\sB^a(E\odot F)$. If $E$ is \hl{full}, that is, if the \hl{range ideal} $\cB_E:=\cls\AB{E,E}$ in $\cB$ is $\cB$, then one may show that this embedding is an isomorphism onto the relative commutant of $\sB^a(E)\odot\id_F$ in $\sB^a(E\odot F)$.

If $(v,w)\mapsto v\cdot w$ is bilinear or sesquilinear operation, then $VW$ is the set $\CB{v\cdot w\colon v\in V,w\in W}$. We \hl{do not} adopt the convention that $VW=\ls VW$ or even $VW=\cls VW$.

%\newpage

\section{The product system associated with an $E_0$--semigroup}\label{psaE0SEC}

Let $\bS$ be one of the (additive) semigroups $\R_+$ or $\N_0$ (with identity $0$). We will refer to $\bS=\R_+$ also as the \hl{continuous time} case and to $\bS=\N_0$ as the \hl{discrete} case. We are mainly interested in the continuous time case. In associating with an \nbd{E_0}semigroup a product system, there is no difference between the discrete and the continuous time case. But knowing how to deal with the discrete case will play a crucial role in showing the converse statement in Sections \ref{AsE0SEC} and \ref{psE0SEC}. For the forward direction in this section, we will discuss first the Hilbert space case and then gradually pass to modules.

Let $\vt=\bfam{\vt_t}_{t\in\bS}$ be an \nbd{E_0}semigroup on the algebra $\sB(H)$ of all adjointable operators on a Hilbert space $H$. Recall that an \hl{\nbd{E_0}semigroup} $\vt$ on a unital \nbd{*}algebra is a semigroup of unital endomorphisms. If the \nbd{*}algebra is $\sB(H)$, then we will require that these endomorphisms are normal, while for the time being we do note pose continuity conditions regarding time dependence of $\vt_t$. We mentioned already, that there are essentially two ways to associate with $\vt$ a product system of Hilbert spaces (\hl{Arveson system}, for short). The first one is Arveson's original construction from \cite{Arv89}, the second one is due to Bhat \cite{Bha96}. However, only the second construction due to Bhat allows for a direct generalization to Hilbert modules. Arveson's construction works nicely only for von Neumann modules and results in a different product system, the \it{commutant system}. Even for Hilbert spaces the results of the constructions need not coincide; see Footnote \ref{ncfn}. Moreover, he results are related to the original \nbd{E_0}semigroup in different ways, namely, one (Arveson) by what we will call a \it{right dilation} and the other (Bhat) by what we will call a \it{left dilation}. Starting with this section we will discuss product systems and their relations with \nbd{E_0}semigroups in terms that correspond rather to Bhat's construction. The generalization of Arveson's approach requires the \it{commutant} of a von Neumann correspondence. We will discuss these things starting from Section \ref{E0AsSEC}.

In \cite{Bha96} Bhat chooses a unit vector $\xi\in H$ and defines the Hilbert subspaces
\beq{\label{Etdef}
E^B_t
~:=~
\vt_t(\xi\xi^*)H
}\eeq
of $H$. (Once for all, for an element $x$ in a space with an inner product, we define the map $x^*\colon y\mapsto\AB{x,y}$. Consequently, $xy^*$ is the \hl{rank-one operator} $z\mapsto x\AB{y,z}$.) It is easy to show that the bilinear maps
\beq{\label{vudef}
(x,y_t)
~\longmapsto~
xy_t
~:=~
\vt_t(x\xi^*)y_t
\text{~~~~~~~~~~~~and~~~~~~~~~~~~}
(x_s,y_t)
~\longmapsto~
x_sy_t
~:=~
\vt_t(x_s\xi^*)y_t
}\eeq
define isometries $v_t\colon H\otimes E^B_t\rightarrow H$ and $u_{s,t}\colon E^B_s\otimes E^B_t\rightarrow E^B_{s+t}$. (Exercise: Check that $u_{s,t}$ is into $E_{s+t}$.) Using a bounded approximate unit of finite-rank operators and normality, one may show that $v_t$ is surjective. (We discuss this in a minute in the more general context; see Equation \eqref{nondeg} and the exercise suggested there.) Now $u_{s,t}$ is just the restriction of $v_t$ to the subspace $E^B_s\otimes E^B_t$ of $H\otimes E^B_t$ and $v_t^*$ maps the subspace $E_{s+t}$ of $H$ into $E_s\otimes E_t$. (Exercise!) This shows that $u_{s,t}$ is onto $E_{s+t}$. We find
\beqn{
(xy_s)z_t
~=~
\vt_t((xy_s)\xi^*)z_t
~=~
\vt_t(\vt_s(x\xi^*)y_s\xi^*)z_t
~=~
\vt_{s+t}(x\xi^*)\vt_t(y_s\xi^*)z_t
~=~
\vt_{s+t}(x\xi^*)(y_sz_t)
~=~
x(y_sz_t)
}\eeqn
and, by restriction, $(x_ry_s)z_t=x_r(y_sz_t)$. Therefore, the family ${E^B}^\otimes=\bfam{E^B_t}_{t\in\bS}$ is an (algebraic) Arveson system\footnote{``Algebraic'' refers to that we are not posing any continuity or measurability condition on ${E^B}^\otimes$.} and the $v_t$ iterate associatively with that product system structure. We call ${E^B}^\otimes$ the \hl{Bhat system associated with $\vt$}.\footnote{Of course, the construction of ${E^B}^\otimes$ depends on $\xi$. But we explain in Proposition \ref{unildprop} that all Arveson systems we obtain from different choices are isomorphic. Moreover, we will single out the result of yet another construction as \hl{the} Bhat system of $\vt$. (That construction has the advantage that it works with choosing a distinguished unit vector. But, even in the Hilbert space case, its simple proof cannot be understood without knowing Hilbert modules; see Remark \ref{absrem}.) If we want to emphasize the unit vector $\xi$, we will say the Bhat system of $\vt$ based on $\xi$.}

In general, whenever for an Arveson system $E^\otimes$ we have a Hilbert space $L\ne\zero$ and a family $w^\otimes$ of unitaries $v_t\colon L\otimes E_t\rightarrow L$ that iterates associatively with the product system structure, we call the pair $(v^\otimes,L)$ a \hl{left dilation} $v^\otimes$ of $E^\otimes$ to $L$.\footnote{\label{canwfn}Note that by associativity and the requirement that $u_{0,0}$ is the canonical identification, it follows that also $v_0$ is the canonical identification. Indeed, suppose $u$ is the unique unitary in $\sB(L)$ such that $v_0(x\otimes 1)=ux$. Then  $ux=v_0(x\otimes 1)=v_0((v_0(u^*x\otimes1))\otimes1)=v_0(u^*x\otimes u_{0,0}(1\otimes 1))=v_0(u^*x\otimes1)=uu^*x=x$, so that $u=\id_L$.} In that case, by setting $\vt_t^v(a):=v_t(a\otimes\id_t)v_t^*$ we define an \nbd{E_0}semigroup $\vt^v$ on $\sB(L)$. (The semigroup property corresponds exactly to the associativity condition.) Moreover, it is easy to check that the Bhat system of $\vt^v$ is $E^\otimes$ by identifying $x_t\in E_t$ with $v_t(\xi\otimes x_t)\in\vt^v_t(\xi\xi^*)L$. (Exercise: Verify that this identification does not only preserve the spaces but also the product system structure.) In the case of the Bhat system ${E^B}^\otimes$ of an \nbd{E_0}semigroup $\vt$ on $\sB(H)$ and the left dilation $v_t$ of ${E^B}^\otimes$ to $H$ as constructed before, it follows from
\beqn{
v_t(a\otimes\id_t)v_t^*(xy_t)
~=~
v_t(a\otimes\id_t)(x\otimes y_t)
~=~
v_t(ax\otimes y_t)
~=~
\vt_t(ax\xi^*)y_t
~=~
\vt_t(a)\vt_t(x\xi^*)y_t
~=~
\vt_t(a)(xy_t)
}\eeqn
that $\vt^v_t=\vt_t$. We summarize:

\bprop\label{Bldprop}
Let $E^\otimes$ be an (algebraic) Arveson system. The problem of finding an \nbd{E_0}semi\-group that has $E^\otimes$ as associated Bhat system is equivalent to the problem of finding a left dilation of $E^\otimes$.
\eprop

Now suppose that $\vt$ is an \nbd{E_0}semigroup acting on $\sB^a(E)$ where $E$ is a Hilbert \nbd{\cB}module. In order to obtain a representation theory of $\sB^a(E)$ on $E$ in analogy with that of $\sB(H)$, we need a condition that replaces normality. The crucial point is that a normal representation of $\sB(H)$ is determined completely by what it does to the rank-one operators. In particular, if the representation is nondegenerate, then already the action of the rank-one operators alone has to be nondegenerate. (For a unital representation of $\sB(H)$, this nondegeneracy condition is equivalent to normality!) We will require that all unital endomorphisms $\vt_t$ of $\sB^a(E)$ are nondegenerate in that sense, that is, we require that for all $t\in\bS$ the set $\vt_t(EE^*)E$ is total in $E$. It can be shown that this is equivalent to say that the unital representation $\vt_t$ is \hl{strict}; see, for instance, \cite{MSS06}.

To begin with, suppose that $E$ has a \hl{unit vector} $\xi$, that is, $\AB{\xi,\xi}=\U\in\cB$. That means, in particular, that $\cB$ is unital and that $E$ is full. We showed in Skeide \cite{Ske02} that, in this case, the whole construction of a product system \it{à la} Bhat \it{cum grano salis} goes through, as before. As in \eqref{Etdef}, we define Hilbert \nbd{\cB}submodules $E_t:=\vt_t(\xi\xi^*)E$ of $E$. The \it{grano salis} we had to add in \cite{Ske02} is the definition of a left action of $\cB$ on $E_t$ that turns it into a correspondence over $\cB$. This left action is
\beqn{
bx_t
~:=~
\vt_t(\xi b\xi^*)x_t.
}\eeqn
(Exercise: Check that this defines a unital representation of $\cB$ by operators on $E_t$.) Once more, the (balanced \nbd{\C}bilinear) mappings in \eqref{vudef} define isometries $v_t\colon E\odot E_t$ and $u_{s,t}\colon E_s\odot E_t\rightarrow E_{s+t}$. (We invite the reader to check that these maps, indeed, preserve inner products.) To see that $v_t$ is surjective, simply observe that the elements of the total subset $\vt_t(EE^*)E$ can be written as
\beq{\label{nondeg}
\vt_t(xy^*)z
~=~
\vt_t(x\xi^*\xi y^*)z
~=~
\vt_t(x\xi^*)\vt_t(\xi y^*)z
~=~
v_t(x\odot\vt_t(\xi y^*)z)
}\eeq
where, clearly, $\vt_t(\xi y^*)z\in E_t$. (Exercise: Go back to the Hilbert space case and give a formal proof of surjectivity under the apparently weaker assumption of normality, modifying the preceding argument suitably.) Of course, also here $\vt^v_t(a):=v_t(a\odot\id_t)v_t^*$ gives back $\vt_t(a)$. Surjectivity of $u_{s,t}$ can be checked as in the Hilbert space case. And by
\beqn{
bu_{s,t}(x_s\odot y_t)
~=~
\vt_{s+t}(\xi b\xi^*)\vt_t(x_s\xi^*)y_t
~=~
\vt_t(\vt_s(\xi b\xi^*)x_s\xi^*)y_t
~=~
\vt_t(bx_s\xi^*)y_t
~=~
u_{s,t}(bx_s\odot y_t)
}\eeqn
we see that the unitaries $u_{s,t}$ are even bilinear.

We summarize: The family $E^\odot=\bfam{E_t}_{t\in\bS}$ with the unitaries $u_{s,t}\in\sB^{a,bil}(E_s\odot E_t,E_{s+t})$ is an (algebraic) \hl{product system} of correspondences over $\cB$. That means, the product $(x_s,y_t)\mapsto x_sy_t:=u_{s,t}(x_s\odot y_t)$ is associative, $E_0=\cB$ and $u_{t,0}$ and $u_{0,t}$ are the canonical identifications. Moreover, the product system is \hl{full} in the sense that each $E_t$ is full, and the pair $(v^\odot,E)$ with $v^\odot=\bfam{v_t}_{t\in\bS}$ is a \hl{left dilation} of $E^\odot$ to $E$. By this we mean that the unitaries $v_t\in\sB^a(E\odot E_t,E)$ iterate associatively with the product system structure \hl{and} that $E$ is full. The \nbd{E_0}semigroup $\vt^v=\bfam{\vt^v_t}_{t\in\bS}$ is the $\vt$ we started with.

\brem
Note that if $(v^\odot,L)$ is a left dilation, then $\cB_L\subset\cB_{E_t}$ so that full $L$ implies that every $E_t$ is full. The condition that $L$ be full replaces the condition $L\ne\zero$ of nontriviality in the Hilbert space case. In fact, the only Hilbert space that is not a full Hilbert \nbd{\C}module is $\zero$. For nonfull $E^\odot$ the concept of left dilation is not defined.
\erem

The idea of left dilation is that, if a product system $E^\odot$ gives rise to an \nbd{E_0}semigroup $\vt^v$ via a left dilation $(v^\odot,L)$, then the \nbd{E_0}semigroup should determine that product system uniquely. By this we mean, if we have another product system with a left dilation to the same $L$ such that the induced \nbd{E_0}semigroups coincide, then the two product systems should be isomorphic. For full $L$ this is a special case of Proposition \ref{unildprop} below. If we would weaken to not necessarily full $L$, then uniqueness fails as soon as $L$ is not full. (An extrem example would be $L=\zero$ to which every product system could be ``dilated''.) If we have a pair $(v^\odot,L)$ that fulfills all conditions of a left dilation except fullness of $L$, then we speak of a left \hl{quasi dilation}. Also a quasi dilation defines an \nbd{E_0}semigroup $\vt^v$ on $\sB^a(L)$.

Every product system $E^\odot$ of correspondences over $\cB$ with a quasi dilation $v^\odot$ to $L$ has a subsystem $F^\odot$ of full correspondences
\beqn{
F_t
~:=~
\bigcap_{t_1+\ldots+t_n=t}\cls\bfam{\cB_LE_{t_n}\cB_L\ldots\cB_LE_{t_1}\cB_L}
}\eeqn
over $\cB_L$. It is easy to check (exercise!) that the restriction of the quasi dilation of $E^\odot$ to that subsystem $F^\odot$ is, now, a left dilation of $F^\odot$ to the full Hilbert \nbd{\cB_L}module $L$ inducing the same \nbd{E_0}semigroup on $\sB^a(L)$. By Proposition \ref{unildprop}, which holds also for nonunital $\cB_L$, such a product system is determined uniquely by $\vt^v$.

We owe the reader to say a few words about the construction of \hl{the} unique product system of an \nbd{E_0}semigroup in the general case. (The reader who is satisfied considering only the full unital case, may skip this and pass to Proposition \ref{unildprop}, immediately.) Again this is nothing but the theory of (strict) representations of $\sB^a(E)$, now in its most general form. The theory of unital strict representations $\vt$ of $\sB^a(E)$ on another Hilbert module $F$ over a possibly different \nbd{C^*}algebra $\cC$ \hl{and} the theory of arbitrary representations on a von Neumann module have been settled in Muhly, Skeide and Solel \cite{MSS06}. In the strict and unital case there is a correspondence $F_\vt$ from $\cB$ to $\cC$ such that $F\cong E\odot F_\vt$ and $\vt(a)$ is just amplification $a\odot\id_{F_\vt}$. In the not necessarily strict case, the representation on a von Neumann module decomposes into a strict unital part, and a part that annihilates the \hl{algebra of finite-rank operators} $\sF(E):=\ls EE^*$ and, therefore, also the \hl{\nbd{C^*}algebra of compact operators} $\cK(E)=\ol{\sF(E)}$.

We repeat briefly what the construction asserts in the case of an \nbd{E_0}semigroup $\vt$ on $\sB^a(E)$ as discussed in Skeide \cite{Ske04p}. To begin with, we do not assume that the Hilbert \nbd{\cB}module $E$ is full. For every $t\in \bS$ we turn $E$ into a correspondence $_tE$ from $\sB^a(E)$ to $\cB$ by defining the left action $ax=\vt_t(a)x$. By the nondegeneracy condition we posed on $\vt_t$, already the action of $\sF(E)\subset\sB^a(E)$ alone on $_tE$ is nondegenerate. In other words, we may also view $_tE$ as correspondence from $\sK(E)$ to $\cB$. We turn $E^*=\CB{x^*\colon x\in E}$ into a correspondence from $\cB$ to $\sB^a(E)$ by defining the inner product $\AB{x^*,y^*}:=xy^*$ and the bimodule action $bx^*a:=(a^*xb^*)^*$. As $\sB^a(E)_{E^*}=\sK(E)$ and $\cB_EE^*$ is total in $E^*$ we may view $E^*$ also as a full correspondence from $\cB_E$ to $\sK(E)$.

It is easy to verify that $E\odot E^*\cong\sK(E)$ via $x\odot y^*\mapsto xy^*$ and $E^*\odot E\cong\cB_E$ via $x^*\odot y\mapsto\AB{x,y}$, as correspondences over $\sK(E)$ and over $\cB_E$, respectively.\footnote{\label{Moritafn}Effectively, as \nbd{\sK(E)}\nbd{\cB_E}correspondence, $E$ is a Morita equivalence from $\sK(E)$ to $\cB_E$ and $E^*$ its \it{inverse} under tensor product. In general, what we nowadays call a \hl{Morita equivalence} from $\cA$ to $\cB$, is a full correspondence $F$ from $\cA$ to $\cB$ for which the canonical homomorphism $\cA\rightarrow\sB^a(F)$ defines an isomorphism onto $\sK(F)$. (Rieffel \cite{Rie74}, who introduced the concept, called $F$ an \it{imprimitivity bimodule}.) With this isomorphism the \nbd{\cB}\nbd{\sK(F)}correspondence $F^*$ can be viewed as \nbd{\cB}\nbd{\cA}correspondence. As $\cB$ is canonically isomorphic to $\sK(F^*)$, also $F^*$ is a Morita equivalence. Almost all what follows, essentially noting that tensoring with a Morita equivalence may be undone by tensoring with its inverse, was already known to Rieffel. What we added to his \it{imprimitivity theorem} \cite[Theorem 6.23]{Rie74}, essentially the representation theory of $\sF(E)$ on a Hilbert space, is the extension to $\sB^a(E)$ and that the representation space may be a Hilbert module.} If we define the correspondence $E_t:=E^*\odot{_tE}$ over $\cB_E$, then
\beqn{
E\odot E_t
~=~
E\odot(E^*\odot{_tE})
~\cong~
(E\odot E^*)\odot{_tE}
~\cong~
\sK(E)\odot{_tE}
~\cong~
_tE
}\eeqn
via
\beqn{
v_t
\colon
x\odot(y^*\odot_tz)
~\longmapsto~
\vt_t(xy^*)z,
}\eeqn
where we denote the elementary tensor of elements $y^*\in E^*$ and $z\in{_tE}$ as $y^*\odot_tz$. Note that this is an isomorphism of correspondences from $\sB^a(E)$ to $\cB$ so that the canonical action $a\odot\id_t$ of $a\in\sB^a(E)$ on the left-hand side corresponds to the canonical action $\vt_t(a)$ of $a$ on the right-hand side. It is readily verified (exercise!) that
\beqn{
(x^*\odot_sy)\odot(x'^*\odot_t y')
~\longmapsto~
x^*\odot_{s+t}(\vt_t(yx'^*)y')
}\eeqn
defines an (obviously, bilinear) unitary $u_{s,t}\colon E_s\odot E_t\rightarrow E_{s+t}$ and that this product is associative. In other words, $E^\odot=\bfam{E_t}_{t\in\bS}$ is a product system of \nbd{\cB_E}correspondences and $v^\odot=\bfam{v_t}_{t\in\bS}$ is a left dilation of $E^\odot$ to $E$ giving back $\vt$ as $\vt^v$. If we want to have a concise construction that works for all \nbd{E_0}semigroups, then we speak about this $E^\odot$ as \hl{the} product system associated with $\vt$.

\brem\label{absrem}
There is a price to be paid, for that this construction works for all \nbd{E_0}semigroups. The members $E_t=E^*\odot{_tE}$ are abstract tensor products, while in every other construction, also Arveson's for Hilbert spaces, the $E_t$ are subspaces of one fixed Banach space (of $E$ in the construction \it{à la} Bhat with a unit vector and of $\sB(H)$ in Arveson's construction; see Section \ref{E0AsSEC}).

Also, the proof is very elegant and simple. But, unlike the other proofs, even in the Hilbert space case it requires some basic knowledge of Hilbert modules. ($H^*$ is a correspondence from $\C$ to $\sK(H)$ and we have to calculate tensor products with this correspondence.) We recommend as an intriguing exercise to redo the theory of normal representation of $\sB(H)$ along the preceding proof. See \cite[Remark 2.2]{Ske05a} and \cite[Example 1.5]{MSS06}.
\erem

After the preceding discussion of the general nonunital and even nonfull case, we will now concentrate on full product systems. What happens if we have two left dilations $({v^1}^\odot,L^1)$ and $({v^2}^\odot,L^2)$ of a full product system $E^\odot$? In the case of Hilbert spaces Arveson's answer (in terms of left dilation) is, the two \nbd{E_0}semigroups $\vt^{v_1}$ and $\vt^{v2}$ must be cocycle conjugate. However, this statement relies on the fact that Arveson's \nbd{E_0}semigroups all act on a $\sB(H)$ where $H$ is infinite-dimensional and separable. In other words, the Hilbert spaces $L_1$ and $L_2$ have the same dimension and, therefore, are isomorphic. The general case of Hilbert modules, is an (interesting) open problem.\footnote{Left (quasi) dilations of the same product system admit direct sums. We expect that it is possible to develop a decomposition theory for left dilations in terms of smallest building blocks.} However, if the Hilbert modules $L^1$ and $L^2$ are isomorphic, then we have the same result as Arveson.\footnote{\label{psMoritafn}We have even more: Suppose $\sB^a(L^1)$ and $\sB^a(L^2)$ are strictly isomorphic, so that there is a Morita equivalence $M$ such that $L^2=L^1\odot M$; see Footnote \ref{Moritafn}. Then two \nbd{E_0}semigroups $\vt^{v^i}$on $\sB^a(L^i)$ $(i=1,2)$ are cocycle conjugate (in the obvious way), if and only if their product systems ${E^i}^\odot$  are \hl{Morita equivalent} via the same Morita equivalence $M$, that is, there is an isomorphism between the product systems ${E^1}^\odot$ and $M\odot {E^2}^\odot\odot M^*:=\bfam{M\odot E^2_t\odot M^*}_{t\in\bS}$. For that, $L^1$ and $L^2$ need not even be modules over the same \nbd{C^*}algebra. See Skeide \cite{Ske04p} for details.} In Skeide \cite{Ske02} we discussed the case with unit vectors. Here we state a slightly more general result directly in terms of left dilations.

\bprop\label{unildprop}
Let $L$ be a Hilbert \nbd{\cB}module. Then for two left dilations $({v^i}^\odot,L)$ of two full product systems ${E^i}^\odot$ $(i=1,2)$ to $L$ the product systems are \hl{isomorphic} (that is, there is a family $u_t\colon E^1_t\rightarrow E^2_t$ of bilinear unitaries fulfilling $u_{s+t}(x_sy_t)=u_s(x_s)u_t(y_t)$ and $u_0=\id_\cB$), if and only if the \nbd{E_0}semigroups $\vt^{v^1}$ and $\vt^{v^2}$ are \hl{cocycle conjugate} (that is, there is a family $\eu_t\in\sB^a(L)$ of unitaries with $\eu_0=\id_L$ such that $\vt^{v^1}_{s+t}(\eu_{s+t})=\eu_s\vt^{v^1}_s(\eu_t)$ and $\vt^{v^2}_t(a)=\eu_t\vt^{v^1}_t(a)\eu_t^*$).
\eprop

\proof
(Sketch.)
If there is a family $u_t$, then $\eu_t:=v^2_t(\id_L\odot u_t){v^1_t}^*$ fulfills the desired properties. (Exercise!) For the backward direction recall that $L$ is necessarily full. Therefore, $E^i_t=L^*\odot L\odot E^i_t$ as correspondences over $\cB$. Moreover, $L\odot E^i_t$ is isomorphic to $_tL^i$ via $v^i_t$, if we  define $_tL^i$ as $L$ when viewed as correspondence from $\sB^a(L)$ to $\cB$ via $\vt^i_t$. Clearly, $\eu_t\colon L\rightarrow L$ is an isomorphism of correspondences when viewed as mapping $_tL^1\rightarrow{_tL^2}$. In other words, $u_t:=\id_{L^*}\odot{v^2_t}^*\eu_t v^1_t$ is an isomorphism $E^1_t=L^*\odot L\odot E^1_t\rightarrow L^*\odot L\odot E^2_t=E^2_t$. These $u_t$ form an isomorphism of product systems. (Exercise!)\qed

\lf
We summarize: Every \nbd{E_0}semigroup $\vt$ on $\sB^a(E)$ leads to a product system $E^\odot$. If $E$ is full, then so is $E^\odot$ and $E^\odot$ is related to $\vt$ via a left dilation $v^\odot$ to $E$ giving back $\vt$ as $\vt^v$. Every other left dilation of $E^\odot$ to $E$ leads to an \nbd{E_0}semigroup cocycle conjugate to $\vt$ and two \nbd{E_0}semigroups on $\sB^a(E)$ have isomorphic product systems, if and only if they are cocycle conjugate.

If $E$ is not full, then we may still associate with an \nbd{E_0}semigroup $\vt$ on $\sB^a(E)$ a product system $E^\odot$. This product system consists of full correspondences over $\cB_E$. So if we simply restrict to $\cB_E$, then we are in the full situation. There is no problem to consider $E_t$ as a (no longer full) correspondence over $\cB$. However, if we insist in having a product system of correspondences over $\cB$, then we must replace $E_0=\cB_E$ with $\cB$. This causes a sharp discontinuity at $t=0$. Such a product system with $\cB_{E_t}=\cB_E\ne\cB$ for $t>0$ will never be continuous in the sense of Definition \ref{cPSdef} below. Product systems where $\cB_{E_t}$ increases in continuous way to $\cB$ may have interesting left quasi dilations if there is a nontrivial subalgebra $\cC$ of $\cB$ such that $\cB_{E_t}\supset\cC$ for all $t$. But the investigation of quasi dilations, so far, has not yet been tackled.

Even if $\bigcap_{t\in\bS}\cB_{E_t}=\zero$ we obtain interesting structures, if we weaken, in the definition of left quasi dilation, unitarity of $v_t$ to isometry.

\bex\label{nonfex}
Put $E=\cB=C_0(0,\infty)$, define the Hilbert submodules $E_t=C_0(t,\infty)$ of $\cB$, and let $\s_t$ denote the usual \it{right shift}. We turn $E_t$ into a correspondence over $\cB$ by defining the left action $b.x_t:=\s_t(b)x_t$. We leave it as an exercise to check that $v_t(x\odot y_t)=\s_t(x)y_t$ defines an isometry $E\odot E_t\rightarrow E$, and that the restriction $u_{s,t}$ to $E_s\odot E_t$ defines a bilinear unitary onto $E_{s+t}$. Clearly, the $u_{s,t}$ turn $E^\odot=\bfam{E_t}_{t\in\bS}$ into a product system and the $v_t$ iterate associatively with the product system structure. Note that $v_t$ is not adjointable, so it is not possible to define an \hl{\nbd{E}semigroup} (that is, a semigroup of not necessarily unital endomorphism) $\vt^v$ on $\sB^a(E)=C_b(0,\infty)$, the multiplier algebra $M(\cB)$ of $\cB$. (In fact, such a semigroup should have to send the identity of $\sB^a(E)$ to the indicator function $\I_{\RO{t,\infty}}$ which is not in $C_b(0,\infty)$.) But, $\s_t$ does define an \nbd{E}semigroup on $\sK(E)=\cB$, and $E^\odot$ may be considered as the product system of that \nbd{E}semigroup. The difficulty disappears for von Neumann modules; see Bhat and Lindsay \cite{BhLi05} and Skeide \cite{Ske04p,Ske06p3} for the obvious generalizations from \nbd{E_0}semigroups to \nbd{E}semigroups.
\eex

The question whether to every full product system $E^\odot$ there exists a left dilation and, therefore, an \nbd{E_0}semigroup that has $E^\odot$ as associated product system is the subject of the following two sections.

\section{Arveson systems and $E_0$--semigroups}\label{AsE0SEC}

One of the most important results about Arveson systems is that every Arveson system is the Arveson system associated with an \nbd{E_0}semigroup; see Section \ref{E0AsSEC} for the terminology we use here. Therefore we refer to this result as the \it{fundamental theorem} about Arveson systems. By Observation \ref{lrob}, below, this is equivalent to say that every Arveson system is the Bhat system associated with an \nbd{E_0}semigroup as described in Section \ref{psaE0SEC} or, by Proposition \ref{Bldprop}, to say that every Arveson system admits a left dilation.

Arveson showed the fundamental theorem in the last of the four articles \cite{Arv89,Arv90a,Arv89a,Arv90}. After laying the basis in \cite{Arv89}, in \cite{Arv90a} he introduced the \it{spectral \nbd{C^*}algebra} of an Arveson system, that is, essentially, the \nbd{C^*}algebra generated by the representing operators of a sufficiently faithful representation (see below) of the Arveson system. A universal property asserts that representations of the Arveson system are in one-to-one correspondence with representations of the spectral algebra. In \cite{Arv89a} he analyzed that sort of representations, the singular representations, that do not lead to \nbd{E_0}semigroups. (In the discrete case, a singular representation corresponds to the defining representation of a Cuntz algebra on the full Fock space.) In \cite{Arv90} he used the precise knowledge of the singular representations to construct an essential representation to every Arveson system, and this is equivalent to constructing an \nbd{E_0}semigroup (Proposition \ref{rdE0prop}). Independently of proving the fundamental theorem, the spectral algebra and the deep analysis of its representations in \cite{Arv90a,Arv89a,Arv90} is interesting in its own right and has been subject to intense research, among others, by Zacharias and Hirshberg \cite{Zac00,Zac00a,HiZa03,Hir04,Hir05,Hir05a}. We shall not discuss the spectral algebra.

In the meantime, there is a proof due to Liebscher \cite{Lie00p1}, which is similarly involved as Arveson's. (See also Footnote \ref{LiFN}.) Since Skeide \cite{Ske06} we have a proof of the fundamental theorem that fits into a few pages. Shortly after, Arveson \cite{Arv06} presented a proof which frees the construction in \cite{Ske06} from a not actually difficult but quite tedious verification of the associativity condition. In Skeide \cite{Ske06a} it is shown that the result of Arveson's construction \cite{Arv06} is, indeed, unitarily equivalent to (a special case of) the construction in \cite{Ske06}. While Arveson's approach \cite{Arv90a,Arv89a,Arv90} via the spectral algebra definitely is not applicable to the case of Hilbert modules, the proof(s) in \cite{Ske06,Arv06,Ske06a} generalize in a (more or less) straightforward way to Hilbert modules; see \cite{Ske06p5} (\nbd{E_0}semigroups) and \cite{Ske06p6} (essential representations) for Hilbert modules and \cite{Ske07p} (in preparation) for von Neumann modules.

\lf
In this section we discuss the case of Hilbert spaces. The versions we have, so far, for modules we discuss in Section \ref{psE0SEC}. We should like to say that we will describe the construction of a left dilation, because it is this construction which generalizes to Hilbert modules. What Arveson constructs (be it in \cite{Arv89a} or in \cite{Arv06}) is an essential representation or, what is the same, a \it{right} dilation; see the end of this section. For Hilbert spaces these concepts may be easily translated into each other; see Observation \ref{lrob}. (In fact, in \cite{Ske06a} we translated \cite{Ske06} into a right dilation in order to compare with \cite{Arv06}. Here we will proceed the other way round.)

The hard problem is the continuous time case $\bS=\R_+$ in absence of so-called \it{units}. A \hl{unit} for an Arveson system $E^\otimes$ is a family $\xi^\otimes=\bfam{\xi_t}_{t\in\bS}$ of elements $\xi_t\in E_t$ that fulfills
\beq{\label{unidef}
\xi_s\xi_t
~=~
\xi_{s+t}
}\eeq
and $\xi_0=1$. Arveson excludes the \hl{trivial} case where $\xi_t=0$ for all $t>0$. We do not want to exclude it at all as a possibility. Nevertheless, in these notes we shall assume tacitly that a unit is nontrivial.\footnote{If we speak about \it{continuous} units, as almost everywhere in the Hilbert modules case, then nontriviality is automatic. Instead, it is a well-known obstacle in semigroup theory that just measurability is not enough.} We say a unit $\xi^\otimes$ is \hl{unital}, if $\AB{\xi_t,\xi_t}=1$ for all $t\in\bS$.

If $E^\otimes$ has a unital unit $\xi^\otimes$, then already Arveson \cite[appendix]{Arv89} constructed a right dilation by an inductive limit. We discuss here the version for left dilations from \cite{BhSk00,BBLS04} that will work also for Hilbert modules. For every $s,t\in\bS$ we define an isometric embedding $E_t\rightarrow E_{s+t}$ by $x_t\mapsto\xi_sx_t$. The $E_t$ together with these embeddings form an inductive system. We denote by $L$ its inductive limit. The factorization $u_{s,t}\colon E_s\otimes E_t\rightarrow E_{s+t}$ under the limit $s\to\infty$ gives rise to a factorization $v_t\colon L\otimes E_t\rightarrow L$. Associativity of the product system structure $u_{s,t}$ guarantees that the $v_t$ form a left dilation of $E^\otimes$. All the $\xi_t\in E_t$ in the inductive limit appear as the same vector $\xi$ and $v_t(\xi\otimes\xi_t)=\xi$. We leave it as an exercise to show that the Bhat system associated with the \nbd{E_0}semigroup $\vt^v$ via $\xi$ is $E^\otimes$.

The problem is that, in the continuous time case, there are loads of product systems without units.\footnote{Powers \cite{Pow87} showed existence of nonspatial \nbd{E_0}semigroups on $\sB(H)$ by rather indirect means. And nonspatiality is equivalent to that the associated Arveson (or Bhat) system is unitless. The first constructive examples are due to Tsirelson \cite{Tsi00p2}. Bhat and Srinivasan \cite{BhSr05} started a systematic investigation of Tsirelson's probabilistic constructions in a more functional analytic way, and discovered a large class of examples.} However, it is always possible to find a unit for a product system in the discrete case $\bS=\N_0$. Simply take any unit vector $\xi_1\in E_1$ and put $\xi_n=\xi_1^n$. Then $\xi^\otimes=\bfam{\xi_n}_{n\in\N_0}$ is a unital unit for $E^\otimes$. Existence of left dilations for discrete Arveson systems is the starting point of the construction in Skeide \cite{Ske06}.

So let $E^\otimes=\bfam{E_t}_{t\in\R_+}$ be an Arveson system. Suppose we have a left dilation $\breve{v}_n$ of the discrete subsystem $\bfam{E_n}_{n\in\N_0}$ of $E^\otimes$ to $\breve{L}$. (This can be the preceding inductive limit construction based on a unit vector $\xi_1\in E_1$, but it need not.) We try now to ``lift'' this left dilation of the discrete subsystem to a left dilation of the whole system.\footnote{\label{LiFN}This very similar to Riesz' proof of Stone's theorem; see Riesz and Sz.-Nagy \cite{RiNa82}. But, there are also similarities to Liebscher's proof in \cite{Lie00p1}. However, the \nbd{E_0}semigroup Liebscher constructs is \it{pure}, while ours is definitely nonpure.} To that goal we consider the direct integrals $\int_a^bE_\alpha\,d\alpha$ ($0\le a<b\le\infty$). We put $L:=\breve{L}\otimes\int_0^1E_\alpha\,d\alpha$. For $t\in\R_+$ we define $n:=\SB{t}$, the unique integer such that $t-n\in\RO{0,1}$. Then the following identifications
\bal{\notag
L\otimes E_t
&
~=~
\breve{L}\otimes\family{\int_0^1E_\alpha\,d\alpha}\otimes E_t
~\cong~
\breve{L}\otimes\int_t^{1+t}E_\alpha\,d\alpha
\\[2ex]\notag
&
~\cong~
\family{\breve{L}\otimes E_n\otimes\int_{t-n}^1E_\alpha\,d\alpha}
\oplus
\family{\breve{L}\otimes E_{n+1}\otimes\int_0^{t-n}E_\alpha\,d\alpha}
\\[2ex]\label{idea}
&
~\cong~
\family{\breve{L}\otimes\int_{t-n}^1E_\alpha\,d\alpha}
\oplus
\family{\breve{L}\otimes\int_0^{t-n}E_\alpha\,d\alpha}
~=~
L
}\eal
define a unitary $v_t\colon L\otimes E_t\rightarrow L$. In the step from the second line to the third one we have made use of the identifications $\breve{v}_n\colon\breve{L}\otimes E_n\rightarrow\breve{L}$ and $\breve{v}_{n+1}\colon\breve{L}\otimes E_{n+1}\rightarrow\breve{L}$ coming from the dilation of $\bfam{E_n}_{n\in{\N_0}}$. Existence of the dilation of the discrete subsystem means that $\breve{L}$ absorbs every tensor power of $E_1$. Just that how many factors $E_1$ have to be absorbed depends on whether $\alpha+t-n$ is bigger or smaller then $1$.

The identifications in \eqref{idea} suggest operations that act directly on sections $\bfam{\breve{x}\otimes y_\alpha}_{\alpha\in\RO{0,1}}$ and give as result again a section in $\bfam{\breve{L}\otimes E_\alpha}_{\alpha\in\RO{0,1}}$. It is a tedious but straightforward verification (one page in \cite{Ske06}) directly on sections that these identifications iterate associatively. But a word need to be said about the direct integrals, because these do no longer make sense without technical conditions.\footnote{\label{didsfn}The direct integrals make sense immediately, if for $d\alpha$ we choose the counting measure, that is, as direct sums. In this case, the \nbd{E_0}semigroup we obtain from the left dilation has no chance to satisfy any reasonable continuity condition, because the left dilation involves a shift of sections. (Also the time shift on $\ell^2\RO{0,1}$ is not weakly continuous.) Anyway: Every algebraic Arveson system admits a left dilation, though not always a continuous one.} The technical conditions on an Arveson system are such that $\bfam{E_t}_{t>0}$ is isomorphic to $(0,\infty)\times H_0$ as a Borel bundle, where $H_0$ is an infinite-dimensional separable Hilbert space such that $E_t\cong H_0$, and such that the product $(x_t,y_s)\mapsto x_ty_s$ is measurable. We will have to speak more about good choices of technical conditions in the module case. Here it is enough to know that the direct integrals gain a sense as $\int_a^bE_\alpha\,d\alpha=L^2(\RO{a,b},H_0)$. The fundamental theorem, in its precise formulation, involves the statement that under the technical conditions on an Arveson system, the constructed \nbd{E_0}semigroup is strongly continuous in time. The self-contained proof in \cite{Ske06} is done by involving Observation \ref{pairob}. (More precisely, we construct in the same way also a right dilation of $E^\otimes$ to $R$, and show that the unitary group $u_t$ on $L\otimes R$ defined in Observation \ref{pairob} is weakly measurable, hence, by separability, strongly continuous. And, this turns over to the \nbd{E_0}semigroup $\vt$.)

Arveson \cite{Arv06} was able to simplify considerably the proof of associativity in the special case, when the left dilation of the discrete subsystem $\bfam{E_n}_{n\in\N_0}$ is obtained as an inductive limit over the unit $\bfam{\xi_1^n}_{n\in\N_0}$ obtained from a unit vector $\xi_1\in E_1$. Starting with that unit vector, Arveson defines the space of \hl{stable sections}, that is, of all locally square integrable sections $\bfam{y_\alpha}_{\alpha\in\R_+}$ that fulfill
\beq{\label{stable}
\xi_1y_\alpha
~=~
y_{\alpha+1}
}\eeq
for all sufficiently big $\alpha$. It is not difficult to show that for two such sections $y$ and $z$ the expression $\int_a^{a+1}\AB{y_\alpha,z_\alpha}\,d\alpha$ is eventually constant, so that
\beqn{
\AB{y,z}
~:=~
\lim_{a\to\infty}\int_a^{a+1}\AB{y_\alpha,z_\alpha}\,d\alpha
}\eeqn
defines a semiinner product on the space of stable sections. The kernel of this semiinner product consists of those stable sections which are eventually $0$ almost everywhere. The quotient is a Hilbert space. In \cite{Ske06a} we have shown that this Hilbert space is canonically isomorphic to $L$.\footnote{Roughly speaking, in the construction of $\breve{L}\otimes\int_0^1E_\alpha\,d\alpha$ one has to interchange inductive limit and direct integral. In other words, one considers an inductive limit over $E_n\otimes\int_0^1E_\alpha\,d\alpha=\int_n^{n+1}E_\alpha\,d\alpha$. This space corresponds to the subspace of stable sections that satisfy \eqref{stable} for all $\alpha\ge n$.} For a stable section $y$ and an element $x_t\in E_t$ Arveson defines the stable section $yx_t$ by setting
\beqn{
(yx_t)_\alpha
~=~
\begin{cases}
y_{\alpha-t}x_t&\alpha\ge t,
\\
0&\text{else}.
\end{cases}
}\eeqn
Then $y\otimes x_t\mapsto yx_t$ defines an isometry that, in the picture $L$, coincides with $v_t$ (from which also surjectivity is immediate). The advantage of Arveson's approach \cite{Arv06} is that associativity is immediate and that it gives an interpretation of the inductive limit in very concrete terms. The construction in \cite{Ske06} is slightly more general. It starts from the well-known observation that it is easy to obtain a left dilation for the discrete subsystem and the basic idea in \eqref{idea} how to transform that dilation into a dilation of the whole system.

We close this section with some explanations about left and right dilations. Whenever for an Arveson system $E^\otimes$ we have a Hilbert space $R\ne\zero$ and a family $w^\otimes$ of unitaries $w_t\colon E_t\otimes R\rightarrow R$ that iterates associatively with the product system structure, we call the pair $(w^\otimes,R)$ a \hl{right dilation} $w^\otimes$ of $E^\otimes$ to $R$.\footnote{Also here it is automatic that $w_0$ is the canonical identification; cf.\ Footnote \ref{canwfn}.}

A \hl{representation} of an Arveson system $E^\otimes$ on a Hilbert space $R$ is a family of maps $\eta_t\colon E_t\rightarrow\sB(R)$ such that
\baln{
\eta_t(x_t)^*\eta_t(y_t)
&~=~
\AB{x_t,y_t}\id_R,
&
\eta_t(x_t)\eta_s(y_s)
&~=~
\eta_{t+s}(x_ty_s).
}\eeqn
A representation is \hl{nondegenerate} (or essential in \cite{Arv89}), if each $\eta_t$ is nondegenerate, that is, if $\cls\eta_t(E_t)R=R$ for all $t\in\bS$.

If we have a representation, then it is easy to check that $w_t\colon x_t\otimes y\mapsto\eta_t(x_t)y$ defines an isometry $w_t\colon E_t\otimes R\rightarrow R$. These isometries iterate associatively with the product system structure. If the representation is nondegenerate, then the $w_t$ form a right dilation. Conversely, if we have a right dilation $w_t$, then $\eta_t(x_t)\colon y\mapsto w_t(x_t\otimes y)$ defines a nondegenerate representation. Therefore, it is the same to speak about a right dilation or a nondegenerate representation. In \cite{Arv89,Arv90a,Arv89a,Arv90} Arveson showed his fundamental theorem by establishing existence of a representation for every Arveson system $E^\otimes$.

\bob\label{lrob}
For an Arveson system $E^\otimes$ we define the \hl{opposite} Arveson system $E'^\otimes$ as the same family of Hilbert spaces (with the same measurable structure) but opposite multiplication $(x_s,y_t)\mapsto y_tx_s$. We may transform a left dilation $(v^\otimes,L)$ of $E^\otimes$ into a right dilation $(w'^\otimes,L)$ of the opposite system and a right dilation $(w^\otimes,R)$ of $E^\otimes$ into a left dilation $(v'^\otimes,L)$ of the opposite system. We simply have to reverse in all tensor products the order of the factors, that is, we put $w'_t(x_t\otimes y):=v_t(y\otimes x_t)$ and $v'_t(x\otimes y_t):=w_t(y_t\otimes x)$.

For correspondences the operation$(x_s,y_t)\mapsto y_tx_s$ will rarely define an isometry. In general, there is no opposite system for product systems of correspondences. However, for product systems of von Neumann correspondences there is the commutant, and the commutant of an Arveson system coincides with its opposite system; see Example \ref{commex}.
\eob

\bob\label{pairob}
Suppose $(v^\otimes,L)$ and $(w^\otimes,R)$ are a left and a right dilation, respectively, of $E^\otimes$. Then $u_t:=(v_t\otimes\id_R)(\id_L\otimes w_t^*)$ defines a unitary group in $\sB(L\otimes R)$. Moreover, the automorphism group $\alpha_t:=u_t\bullet u_t^*$ leaves invariant $\sB(L)\otimes\id_R$ for $t\ge0$ and $\id_L\otimes\sB(R)$ for $t\le0$. Then the restriction $\alpha_t$ $(t\ge0)$ to $\sB(L)\otimes\id_R\cong\sB(L)$ defines an \nbd{E_0}semigroup $\vt$ on $\sB(L)$ which has $E^\otimes$ as associated Bhat system, while the restriction $\alpha_t$ $(t\le0)$ to $\id_L\otimes\sB(R)\cong\sB(R)$ defines an \nbd{E_0}semigroup $\theta$ on $\sB(R)$ which has $E'^\otimes$ as associated Bhat system. Canceling from the last phrase all statements about product systems, this is exactly the situation when Powers and Robinson \cite{PoRo89} say $\vt$ and $\theta$ are \hl{paired}.

Note that the \nbd{E_0}semigroup $\theta$ on $\sB(R)$ has $E^\otimes$ as associated Arveson system in the sense of \cite{Arv89} and $\vt$ has $E'^\otimes$ as associated Arveson system. We explain this in Section \ref{E0AsSEC}.
\eob

\section{Continuous product systems and $E_0$--semigroups}\label{psE0SEC}

Let $E^\odot=\bfam{E_t}_{t\in\bS}$ be a product system of correspondences over a unital \nbd{C^*}algebra $\cB$. Like in the Hilbert space case, a \hl{unit} for $E^\odot$ is a family $\xi^\odot=\bfam{\xi_t}_{t\in\bS}$ of elements $\xi_t\in E_t$ fulfilling \eqref{unidef} and $\xi_0=\U$. We \hl{do not} define what a unit is, if $\cB$ is nonunital! The unit is \hl{unital}, if $\AB{\xi_t,\xi_t}=\U$ for all $t\in\bS$. The construction of an \nbd{E_0}semigroup from a unital unit works as for Hilbert spaces: We define isometric embeddings $E_t\rightarrow E_{s+t}$ by $x_t\mapsto\xi_sx_t$. Then, the inductive limit $L$ factors as $v_t\colon L\odot E_t\rightarrow L$ via a left dilation $v^\odot$.

\brem\label{nonrrem}
For Hilbert spaces the isometric embedding could be defined as $x_t\mapsto x_t\xi_t$, leading to a right dilation. In fact, this is what Arveson did in the appendix of \cite{Arv89} and what we did in \cite{Ske06a} in the discrete case when we compared Arveson's construction of a right dilation in \cite{Arv06} with \cite{Ske06}. From the beginning, such a construction of a right dilation can not be done for Hilbert modules: $x_t\mapsto x_t\xi_t$ is, in general, not an isometry. More precisely, it is an isometry, if and only if the unital unit $\xi^\odot$ is \hl{central}, that is, if $b\xi_t=\xi_tb$ for all $t\in\bS,b\in\cB$. Product systems that admit a central unital unit are classified as \hl{spatial} (Skeide \cite{Ske06d}) and admit classification results parallel to those for spatial Arveson systems. Product systems of von Neumann modules that admit (continuous) units are spatial automatically (Barreto, Bhat, Liebscher and Skeide \cite{BBLS04}), while for Hilbert modules existence of a continuous unit is not enough to guarantee spatiality. We see, the problem of constructing an \nbd{E_0}semigroup is difficult only for nonspatial product systems.
\erem

If in a discrete product system $E^\odot=\bfam{E_n}_{n\in\N_0}$ the member $E_1$ has a unit vector, then we may construct a left dilation of that discrete product system. However, there are discrete product systems where no member except $E_0$ has a unit vector. It is one of the main results of Skeide \cite{Ske04p} to show that every \it{full} discrete product system of correspondences over a unital \nbd{C^*}algebra admits a left dilation; see Footnote \ref{unifn}.

Now let us discuss the continuous time case. Let $E^\odot=\bfam{E_t}_{t\in\R_+}$ be a full product system of correspondences over a unital \nbd{C^*}algebra $\cB$. By \cite{Ske04p} we may choose a left dilation $(\breve{v}^\odot,\breve{L})$ of the discrete subsystem $\bfam{E_t}_{t\in\N_0}$. We would like to proceed as in Section \ref{AsE0SEC}, defining the direct integrals $\int_a^b E_\alpha\,d\alpha$ and $L=\breve{L}\odot\int_0^1 E_\alpha\,d\alpha$ so that \eqref{idea} had a chance to define a left dilation. Without posing precise technical conditions, this works only for direct sums, leading to noncontinuous \nbd{E_0}semigroups; cf.\ Footnote \ref{didsfn}.

To motivate the technical conditions we will pose, we start from what, in the end, we wish to have. Given a full product system $E^\odot=\bfam{E_t}_{t\in\R_+}$ of correspondences over a unital \nbd{C^*}algebra, we wish to have a full Hilbert module $L$ and a left dilation $v_t\colon L\odot E_t\rightarrow L$ such that the \nbd{E_0}semigroup $\vt^v$ on $\sB^a(L)$ is \hl{strongly continuous}, that is, for all $a\in\sB^a(L)$ and all $x\in L$, the function $t\mapsto\vt^v_t(a)x$ is continuous. (As usual with semigroups, it is sufficient to require strong continuity around $0$. Since the $\vt^v_t$ are bounded uniformly, it is also sufficient to check continuity for $a$ and $x$ from total subsets of $\sB^a(L)$ and $L$, respectively.)

So suppose we have a strongly continuous \nbd{E_0}semigroup $\vt$ on $\sB^a(E)$ where $E$ is a full Hilbert module over a unital \nbd{C^*}algebra $\cB$. What can we say about the associated product system? How is continuity of $\vt$ reflected by the \it{bundle structure} of the product system? As a first step, we have to fix a version of the product system and of the left dilation relating it to $\vt$. We have to take into account that the structures we derive might depend on that choice. To have a start, let us suppose that $E$ has a unit vector $\xi$ and construct product system $E^\odot$ and left dilation $v_t$ of that product system from the unit vector $\xi$ \it{à la} Bhat. The essential observation is that in this approach all $E_t$ are identified as submodules of $E$. Moreover, for every $x\in E$ the function $t\mapsto\vt_t(\xi\xi^*)x\in E_t\subset E$ is continuous. Of course, if $x=y_t\in E_t$, then the section $t\mapsto x_t:=\vt_t(\xi\xi^*)x$ assumes the value $x_t=y_t$ at $t$. It is not difficult to check that, whenever we have two sections $x,y$ of $E^\odot$ such that the functions $t\mapsto x_t\in E_t\subset E$ and $t\mapsto y_t\in E_t\subset E$ are continuous, then also the function $(s,t)\mapsto x_sy_t\in E_{s+t}\subset E$ is continuous; see Skeide \cite{Ske03b}. This motivates the following definition from \cite{Ske03b,Ske06p5}.

\bdefi\label{cPSdef}
Let $E^\odot=\bfam{E_t}_{t\in\R_+}$ be a product system of correspondences over a \nbd{C^*}algebra $\cB$ with a family $i=\bfam{i_t}_{t\in\R_+}$ of isometric embeddings $i_t\colon E_t\rightarrow\wh{E}$ into a Hilbert \nbd{\cB}module $\wh{E}$. Denote by
$$
CS_i(E^\odot)
~=~
\BCB{x=\bfam{x_t}_{t\in\R_+}\colon x_t\in E_t,t\mapsto i_tx_t\mbox{~is continuous}}
$$
the set of \hl{continuous sections} of $E^\odot$ (with respect to $i$). We say $E^\odot$ is \hl{continuous} (with respect to $i$), if the following conditions are satisfied.
\begin{enumerate}
\item\label{csdef}
For every $y_t\in E_t$ we can find a continuous section $x\in CS_i(E^\odot)$ such that $x_t=y_t$.

\item \label{cpsdef}
For every pair $x,y\in CS_i(E^\odot)$ of continuous sections the function
$$
(s,t)
~\longmapsto~
i_{s+t}(x_sy_t)
$$
is continuous.
\end{enumerate}
We say two embeddings $i$ and $i'$ have the same \hl{continuous structure}, if $CS_i(E^\odot)=CS_{i'}(E^\odot)$.
\edefi

\brem\label{Atrivrem}
The definition says, roughly speaking, that $E^\odot$ is a subbundle of the trivial Banach bundle $\RO{0,\infty}\times\wh{E}$. Note that this is even weaker than Arveson's requirement that the part $t>0$ of an Arveson system be \it{Borel isomorphic} to the trivial bundle $(0,\infty)\times H_0$. Of course, this bundle is also a Banach bundle, and the condition just means that the Borel structure of an Arveson system is that induced from the continuous structure of a trivial Banach bundle. We allow even for subbundles.

The only difference between a \it{continuous} Arveson system (that is, an algebraic Arveson system that is continuous in the sense of Definition \ref{cPSdef}) and a \it{measurable} Arveson system (that is, an Arveson system in the sense of Arveson \cite{Arv89}) consists in whether multiplication is required continuous or just measurable. (Note that in \cite{Ske06} we used Arveson's condition, that is, just measurable multiplication.) If, for a general product system, we pose just measurability as assumption, then the construction we describe in the sequel should provide us with a (weakly) measurable \nbd{E_0}semigroup. Only under separability assumptions this will be enough to show strong continuity. We did not yet put into practice a measurable version.
\erem

In Definition \ref{cPSdef} we do not require that $E^\odot$ is full, nor that $\cB$ is unital. However, if $\cB$ is unital, then we have the following lemma from \cite{Ske06}.

\blem\label{cflem}
If $\cB$ is unital, then a continuous product system $E^\odot$ of correspondences over $\cB$ contains a continuous section $\zeta\in CS_i(E^\odot)$ that consists entirely of unit vectors and fulfills $\zeta_0=\U$. In particular, every $E_t$ contains a unit vector (and, therefore, is full).
\elem

The proof relies on the fact that the invertible elements of $\cB$ form an open subset (so that there is a continuous section that consists of unit vectors at least for all sufficiently small $t$), and on the fact that the tensor product of unit vectors is again a unit vector (so that small pieces of that section can be used to compose a global continuous section of unit vectors).

We see that not only the continuous product system is full automatically\footnote{Note that this may fail, if $\cB$ is nonunital. Indeed, the product system from Example \ref{nonfex} with the canonical embedding $C_0(t,\infty)\rightarrow C_0(0,\infty)$ is continuous. But, none of the $E_t$ is full but $E_0$.}, but that all of its members contain a unit vector. That is, we may not only use a unit vector in $E_1$ to construct a left dilation of the discrete subsystem, but we may even adapt Arveson's construction of an \nbd{E_0}semigroup word by word as described in Section \ref{AsE0SEC}.

\brem
Note that, like in the construction of the dilation of the discrete subsystem, also for imitating Arveson's proof it is indispensable that we define the space of \it{stable} sections as in \eqref{stable}, by multiplying the unit vector from the left. Multiplying from the right (as in \cite{Arv06}) will cause that the inner product $\int_a^{a+1}\AB{y_\alpha,z_\alpha}\,d\alpha$ is no longer eventually constant.
\erem

Arveson's proof in \cite{Arv06}, once understood the idea, is easier to carry out than our proof in \cite{Ske06}. This is, why in \cite{Ske06p5} we followed Arveson's road to prove the \it{fundamental theorem} for continuous product systems of correspondences over a unital \nbd{C^*}algebra. The result:

\bthm\label{mthm}
Every continuous product system of correspondences over a unital \nbd{C^*}algebra is the continuous product system associated with a strictly continuous \nbd{E_0}semigroup that acts on the algebra of all adjointable operators on a Hilbert module with a unit vector.
\ethm

Note that the technically most difficult part in the proof of that theorem is to show that the continuous structure induced by the constructed \nbd{E_0}semigroup and the unit vector, is the same we started with.

\brem
So suppose, once more, we have a strongly continuous \nbd{E_0}semigroup $\vt$ on $\sB^a(E)$ where $E$ is a full Hilbert module over a unital \nbd{C^*}algebra $\cB$. If $E$ has a unit vector, then it is not difficult to show that the continuous structure induced on the product system $E^\odot$ associated with $\vt$ does not depend on the choice of that unit vector. As a sort of surprise, if the members of a product system (continuous or not) derived from an \nbd{E_0}semigroup on $\sB^a(E)$ have unit vectors, this does not mean that the full Hilbert module $E$ need have a unit vector.\footnote{As a trivial example, take the Hilbert \nbd{M_2}modules $E=\C_2$ where we put $\C_n:=(\C^n)^*$. The only \nbd{E_0}semigroup $\vt$ on $\sB^a(E)=\C$ is $\vt_t=\id_\C$, its product system simply $E_t=M_2$ with multiplication as identification $E_s\odot E_t=E_{s+t}$. Also the left dilation $v_t$ is simply the canonical identification $E\odot M_2=E$. Nontrivial examples may be obtained by working in the present one via direct sum constructions.} However, if a full $E$ does not have a unit vector, then a little lemma from \cite{Ske04p} asserts that a finite direct sum $E^n$ of copies of $E$ will have a unit vector. (The proof uses, again, that the invertibles form an open subset, and Cauchy-Schwartz inequality: If $E$ is full, then $\U$ is approximated by finite sums $b=\sum_{i=1}^n\AB{x_i,y_i}$. If the approximation is sufficiently good, then $b$ is invertible. $b$ can be interpreted as inner product of elements in $X,Y\in E^n$. A simple application of Cauchy-Schwartz inequality shows that also $\AB{X,X}$ must be invertible so that $X\sqrt{\AB{X,X}^{-1}}$ is a unit vector in $E^n$.)\footnote{\label{unifn}This result is key in the proof of \cite{Ske04p} that every full discrete product system $E^\odot=\bfam{E_n}_{n\in\N_0}$ of correspondences over a unital \nbd{C^*}algebra admits a left dilation. In fact, if $E_1^n$ contains a unit vector, then intuitively also $M_{n\cdot\infty,\infty}(E_1)=M_\infty(E)$ should contain a unit vector. The problem is that $M_\infty(E_1)$, usually, is not big enough, but a suitable \it{strict completion} is. This problem does not appear in the version for von Neumann modules. The price to be paid is that the analogue of the lemma for von Neumann modules requires a direct sum $E_1^\en$ with arbitrary cardinality $\en$. Now $M_\en(E_1)$ is nothing but $\cB^\en\odot E_1\odot\cB_\en$. In other words, the tensor powers of $M_\en(E_1)$ form a product system $\cB^\en\odot E^\odot\odot\cB_\en$ Morita equivalent to $E^\odot$ in the sense of Footnote \ref{psMoritafn}. Once a unit vector in $M_\en(E_1)$ is established, we find a left dilation of $\cB^\en\odot E^\odot\odot\cB_\en$. And one of the major results of \cite{Ske04p} asserts that a product system admits a left dilation, if (and only if) it is Morita equivalent to another product system that admits a left dilation.} The strongly continuous \nbd{E_0}semigroup $\vt$ on $\sB^a(E)$ may be lifted to a strongly continuous \nbd{E_0}semigroup on $\sB^a(E^n)=M_n(\sB^a(E))$ (acting pointwise with $\vt$ on the matrix elements), having the same product system $E^\odot$ as $\vt$. Now $E^\odot$ can be induced from a unit vector $E^n$. Also here the continuous structure depends neither on how big $n$ is, nor on which unit vector we choose.

It is an open problem how to define a continuous structure on $E^\odot$ without reference to a unit vector. A solution might be to give a definition  of continuous product systems in terms of Banach bundles (as Hirshberg's \cite{Hir04} for Borel bundles), that is, by giving explicitly a set of sections that are supposed to be continuous and that determine the structure of the bundle. The product system will, then, be considered as obtained via \cite{MSS06,Ske04p}. A candidate for the generating set of continuous sections would be the set
\beqn{
\BCB{\bfam{x^*\odot_t y}_{t\in\R_+}~|~x,y\in E}.
}\eeqn
It is unclear in how far a definition of continuous product system as Banach bundle (generated by a compatible set of continuous sections) is already sufficient to run through the proof of \cite{Ske06p5}. It might be necessary to find a further condition that substitutes the condition being a subbundle of a trivial bundle.
\erem

\brem
We mentioned in Remark \ref{nonrrem} that Arveson's construction of a right dilation (that is, of an essential representation) for the Hilbert space case, with tensoring a unit vector from the right, fails for modules. However, in Skeide \cite{Ske06p6} we pointed out that the construction, indeed, can be saved if we tensor something different from the right. This ``something different'' is a unit vector not for $E_1$ but for the member $E'_1$ of the commutant of the product system $E^\odot$. And an element of $E_t$ can be tensored with en element from $E'_s$ in a reasonable way. We explain the commutant of a product system in Section \ref{ApsSEC}. But we do not have the space to explain any detail (in particular, how the \nbd{C^*}setting of this section fits into the von Neumann--setting of Section \ref{ApsSEC}) and refer the reader to \cite{Ske06p6}. However, we mention that the existence of a unit vector in $E'_1$ follows by Hirshberg's result \cite{Hir05a} that every full discrete product system with faithful left action admits a right dilation. This result is \it{dual} to \cite{Ske04p} in the sense of commutant. And, as a matter of fact, the condition that the left actions of the product system be faithful is dual to fullness under commutant. For right dilations it is as indispensable as fullness is for left dilations.
\erem

\section{The Arveson system of an $E_0$--semigroup on $\sB(H)$}\label{E0AsSEC}

Preparing Section \ref{ApsSEC}, in this section we review Arveson's construction of an Arveson system from an \nbd{E_0}semigroup $\vt$ on $\sB(H)$ and compare it with Bhat's. Arveson defines for every $\vt_t$ the intertwiner space
\beqn{
E^A_t
~:=~
\bCB{x'_t\in\sB(H)\colon\vt_t(a)x'_t=x'_ta~(a\in\sB(H))}.
}\eeqn
It is readily verified that $x'^*_ty'_t\in\sB(H)'=\C\id_H$, so that $x'^*_ty'_t=\AB{x'_t,y'_t}\id_H$ defines a an inner product turning $E^A_t$ into a Hilbert space. Similarly, it is an easy exercise to check that the bilinear mappings
\baln{
(x'_t,y)
&
~\longmapsto~
x'_ty
&
(x'_t,y'_s)
&
~\longmapsto~
x'_ty'_s
}\ealn
define isometries $w'_t\colon E^A_t\otimes H\rightarrow H$ and $u'_{t,s}\colon E^A_t\otimes E^A_s\rightarrow E^A_{t+s}$. Less obvious is that $w'_t$ is surjective. (Here, normality of $\vt_t$ is essential; see Lemma \ref{MuSolem} in the more general context.) But, from surjectivity of $w'_t$, $w'_s$ and $w'_{t+s}$ it is immediate (exercise!) that also $u'_{t,s}$ is surjective. From $(x'_ty'_s)z'_r=x'_t(y'_sz'_r)$ we see that the $u'_{t,s}$ iterate associatively on multiple tensor products. In other words, the family ${E^A}^\otimes=\bfam{E^A_t}_{t\in\bS}$ forms an (algebraic) Arveson system. We call ${E^A}^\otimes$ the \hl{Arveson system associated with $\vt$}. But, we have more than just the Arveson system of $\vt$. From $(x'_ty'_s)z=x'_t(y'_sz)$ it follows that the $w'_t$ define a right dilation $(w'^\otimes,H)$ of ${E^A}^\otimes$.

In general, if $(w^\otimes,R)$ is a right dilation of $E^\otimes$, by setting $\vt_t^w(a):=w_t(\id_t\otimes a)w_t^*$ we define an \nbd{E_0}semigroup $\vt^w$ on $\sB(R)$. It is easy to check that the Arveson system of $\vt^w$ is $E^\otimes$ by identifying $x_t\in E_t$ with the intertwiner
\beqn{
w_t(x_t\otimes\id_R)
\colon
z
~\longmapsto~
w_t(x_t\otimes z).
}\eeqn
In the case of the Arveson system ${E^A}^\otimes$ of an \nbd{E_0}semigroup $\vt$ on $\sB(H)$ and the right dilation $w'_t$ of ${E^A}^\otimes$ to $H$ as constructed before, it follows from
\beqn{
(w'_t(\id_t\otimes a)w'^*_t)(x'_ty)
~=~
w'_t(\id_t\otimes a)(x'_t\otimes y)
~=~
w'_t(x'_t\otimes ay)
~=~
x'_tay
~=~
\vt_t(a)(x'_ty)
}\eeqn
that $\vt^{w'}_t=\vt_t$. We summarize:

\bprop\label{rdE0prop}
Let $E^\otimes$ be an (algebraic) Arveson system. The problem of finding an \nbd{E_0}semi\-group that has $E^\otimes$ as associated Arveson system is equivalent to the problem of finding a right dilation of $E^\otimes$ or, equivalently, a nondegenerate representation.
\eprop

Let us return to the \nbd{E_0}semigroup $\vt$ on $\sB(H)$ and ask what the relation is between ${E^A}^\otimes$ and ${E^B}^\otimes$. Of course, the dimension of the multiplicity space of an endomorphism $\vt_t$ is unique, no matter whether we factor it out to the left (right dilation) or to the right (left dilation). Therefore, $E^A_t\cong E^B_t$ as Hilbert spaces. But what happens to the identifications $u'_{t,s}$ and $u_{s,t}$? Can we identify $E^A_t$ with $E^B_t$ in such a way that these identifications are preserved? At this point the reader will have noticed that in the discussion of the Arveson system of $\vt$ we discussed an identification of $E^A_t\otimes E^A_s$ with $E^A_{t+s}$, while for ${E^B}^\otimes$ we chose the opposite direction of times. From this, the reader might guess what the answer will be: ${E^A}^\otimes$ and ${E^B}^\otimes$ turn out to be anti-isomorphic, that is , ${E^A}^\otimes$ is the opposite product system of ${E^B}^\otimes$. The two need not be isomorphic.\footnote{\label{ncfn}Tsirelson \cite{Tsi00p1} provided us with an explicit example of a product system that is not isomorphic to its opposite product system.}

Let us see why ${E^A}^\otimes$ and ${E^B}^\otimes$ are anti-isomorphic. Let $x_t\in E^B_t$ and define $x'_t\colon h\mapsto\vt_t(h\xi^*)x_t$ $(=v_t(h\otimes x_t))$. It follows that
\beqn{
\vt_t(a)x'_th
~=~
\vt_t(a)\vt_t(h\xi^*)x_t
~=~
\vt_t(ah\xi^*)x_t
~=~
x'_tah,
}\eeqn
so that $x'_t\in E^A_t$. Clearly, $u_t\colon x_t\mapsto x'_t$ is an isometry. Now let $x'_t\in E^A_t$ and define $x_t=x'_t\xi$. Then
\beqn{
\vt_t(\xi\xi^*)x_t
~=~
\vt_t(\xi\xi^*)x'_t\xi
~=~
x'_t(\xi\xi^*)\xi
~=~
x'_t\xi
~=~
x_t,
}\eeqn
so that $x_t\in E^B_t$. Clearly, $u'_t\colon x'_t\mapsto x_t$ is an isometry. Moreover,
\beqn{
(u_tu'_t x'_t)h
~=~
(u_tx'_t\xi)h
~=~
\vt_t(h\xi^*)x'_t\xi
~=~
x'_t(h\xi^*)\xi
~=~
x'_th.
}\eeqn
In other words, $u_t$ and $u'_t$ are a pair of inverse unitaries. Now let us see what $u'_{t+s}$ does to a tensor product $u'_{t,s}(x'_t\otimes y'_s)$.
\bmun{
u'_{t+s}u'_{t,s}(x'_t\otimes y'_s)
~=~
u'_{t+s}(x'_ty'_s)
~=~
(x'_ty'_s)\xi
~=~
x'_t(y'_s\xi)
~=~
x'_t(u'_sy'_s)
~=~
x'_t((u'_sy'_s)\xi^*)\xi
\\
~=~
\vt_t((u'_sy'_s)\xi^*)x'_t\xi
~=~
\vt_t((u'_sy'_s)\xi^*)(u'_tx'_t)
~=~
u_{s,t}((u'_sy'_s)\otimes(u'_tx'_t)).
}\emun
In other words, $u'_{t+s}u'_{t,s}=u_{s,t}(u'_s\otimes u'_t)$, that is, the family $u'_t$ establishes an anti-isomorphism ${E^A}^\otimes\rightarrow{E^B}^\otimes$.

\section{$E_0$--Semigroups and product systems \it{à la} Arveson: Commutants of von Neumann correspondences}\label{ApsSEC}

Of course, the construction in Section \ref{psaE0SEC} of the product system associated with an \nbd{E_0}semi\-group on $\sB^a(E)$ works also if $E$ is a von Neumann module. (After all a von Neumann module is also a Hilbert module.) In presence of a unit vector it is even clear that the product system consists of von Neumann modules. (The ranges of projections on von Neumann modules are von Neumann modules.) The point is that in the assumptions on the \nbd{E_0}semigroup it is sufficient that the endomorphisms of the von Neumann algebra $\sB^a(E)$ be only normal, not necessarily strict. In Skeide \cite{Ske05a} we provided a generalization of Bhat's approach (without unit vectors, not along the lines of \cite{MSS06}) that works for every von Neumann module (and, of course, gives a product system isomorphic to that constructed along the lines of \cite{MSS06}). As all modifications to be done are plain, we do not discuss them here.

The approach we want to discuss here, is the generalization of Arveson's approach, as discovered in Skeide \cite{Ske03c} together with the commutant of von Neumann correspondences. To that goal we have to spend some time to review the necessary notions and facts about von Neumann modules, von Neumann correspondences and their commutants. The correspondence between a von Neumann algebra and its commutant is bijective. In order that this desirable property remains true for commutants of von Neumann correspondences and does not degenerate to an equivalence, we have to choose our categories carefully. The correct category that allows to view the commutant as a bijective functor is the category of \it{concrete} von Neumann correspondences; Skeide \cite{Ske06b}. In the sequel, we discuss only the case relevant to us, namely, correspondences over $\cB$. (See also Remark \ref{gencommrem}.)

\lf
Before we can speak about concrete von Neumann correspondences, we have to speak about concrete von Neumann modules. Recall that a von Neumann algebra is a strongly closed \nbd{*}algebra $\cB\subset\sB(G)$ of operators acting nondegenerately on a Hilbert space $G$. As usual, by $\cB'\subset\sB(G)$ we denote the commutant of $\cB$. Similarly, a \hl{concrete von Neumann \nbd{\cB}module} is a subset $E$ of $\sB(G,H)$, where $H$ is another Hilbert space, such that
\begin{enumerate}
\item\label{cvnM1}
$E$ is a right \nbd{\cB}submodule of $\sB(G,H)$, that is, $E\cB\subset E$,

\item\label{cvnM2}
$E$ is a pre-Hilbert \nbd{\cB}module with inner product $\AB{x,y}=x^*y$, that is, $E^*E\subset\cB$,

\item\label{cvnM3}
$E$ acts nondegenerately on $G$, that is, $\cls EG=H$, and

\item\label{cvnM4}
$E$ is strongly closed in $\sB(G,H)$.
\end{enumerate}
If we wish to underline the Hilbert space $H$, we will also write the pair $(E,H)$ for the concrete von Neumann \nbd{\cB}module.  One may show (see Skeide \cite{Ske00b,Ske05c}) that a subset $E$ of $\sB(G,H)$ fulfilling \ref{cvnM1}--\ref{cvnM3} (that is, $E$ is a concrete pre-Hilbert \nbd{\cB}module) is a concrete von Neumann \nbd{\cB}module, if and only if $E$ is self-dual\footnote{Recall that a Hilbert \nbd{\cB}module is \hl{self-dual}, if every bounded right linear map $E\rightarrow\cB$ has the form $x\mapsto\AB{y,x}$ for a suitable $y\in E$.}, that is, if and only if $E$ is a \nbd{W^*}module over the von Neumann algebra $\cB\subset\sB(G)$ considered as a \nbd{W^*}algebra. By $\cvN_\cB$ we denote the \hl{category of concrete von Neumann \nbd{\cB}modules} with the adjointable maps $a\in\sB^a(E_1,E_2)$ as morphisms. The definition of \it{concrete} von Neumann modules and their category is due to Skeide \cite{Ske06b}, while the definition of von Neumann modules is due to Skeide \cite{Ske00b}; see Footnote \ref{vNfn}.

Identifying $xg\in H$ with $x\odot g\in E\odot G$, we see from \ref{cvnM3} that $H$ and $E\odot G$ are canonically isomorphic.\footnote{\label{vNfn}In fact, if $E$ is a pre-Hilbert module over a pre-\nbd{C^*}algebra $\cB\subset\sB(G)$, then one may construct the Hilbert space $E\odot G$ with an embedding $x\mapsto L_x\in\sB(G,E\odot G)$ where we put $L_xg:=x\odot g$, transforming $E$ into a concrete pre-Hilbert \nbd{\cB}module $(E,E\odot G)$. For a von Neumann algebra $\cB\subset\sB(G)$ we defined in \cite{Ske00b} that $E$ is a von Neumann \nbd{\cB}module, if its image in $\sB(G,E\odot G)$ is strongly closed. Of course, in that way also a \nbd{W^*}module over a \nbd{W^*}algebra $M$ may be turned into a von Neumann module after choosing a faithful normal unital representation of $M$ on a Hilbert space $G$, thus, turning $M$ into a von Neumann algebra.} Giving $E$ as a subset of $\sB(G,H)$ from the beginning, is crucial for that the commutant, later on, will be bijective. The fact that $H$ is canonically isomorphic to the tensor product $E\odot G$ is, however, by far more inspiring from the algebraic point of view.

For instance, every adjointable operator $a\in\sB^a(E_1,E_2)$ amplifies to an operator $a\odot\id_G\in\sB(E_1\odot G,E_2\odot G)$. Consequently, $a$ gives rise to and is determined uniquely by an operator in $\sB(H_1,H_2)$ that acts as $x_1g\mapsto(ax_1)g$. We shall denote this operator by the same symbol $a$ and identify in that way $\sB^a(E_1,E_2)$ as a subset of $\sB(H_1,H_2)$. It is easy to show that $\sB^a(E_1,E_2)$ is strongly closed in $\sB(H_1,H_2)$. In particular, $\sB^a(E)\subset\sB(H)$ is a von Neumann algebra acting on $H$.

Those operators on the second factor $G$ in $E\odot G$ that embed into $\sB(E\odot G)$ are the \nbd{\cB}\nbd{\C}linear operators on $G$. Of course, $\sB^{bil}(G)=\cB'$ is nothing but the commutant of $\cB$. So, the (clearly, normal and nondegenerate) representation $b'\mapsto\id_E\odot b'$ of $\cB'$ on $E\odot G$ gives rise to a normal nondegenerate representation $\rho'$ of $\cB'$ on $H$ which acts as $\rho'(b')xg=xb'g$. We call $\rho'$ the \hl{commutant lifting} associated with $E$.

From the commutant lifting $\rho'$ we obtain back $E$ as the space
\beq{\label{CB'def}
E
~=~
C_{\cB'}(\sB(G,H))
~:=~
\bCB{x\in\sB(G,H)\colon\rho'(b')x=xb'~(b'\in\cB')}
}\eeq
of \it{intertwiners} for the natural actions of $\cB'$. This was known already to Rieffel \cite{Rie74a}. In \cite{Ske05c} we proved it by simply calculating the double commutant of the linking von Neumann algebra in $\sB(G\oplus H)$:
\beqn{
\fMatrix{\cB&E^*\\E&\sB^a(E)}''
~=~
\left\{\fMatrix{b'&0\\0&\rho'(b')}\colon b'\in\cB'\right\}'
~=~
\fMatrix{\cB&C_{\cB'}(\sB(H,G))\\C_{\cB'}(\sB(G,H))&\rho'(\cB')'}.
}\eeqn
This proof also shows that the commutant $\rho'(\cB')'$ of the range of $\rho'$ in $\sB(H)$ may be identified with the von Neumann algebra $\sB^a(E)\subset\sB(H)$. By doing the computation for $E=E_1\oplus E_2$ one also shows that $\sB^a(E_1,E_2)$ is just $\sB^{bil}(H_1,H_2)$, the space of operators that \it{intertwine} the commutant liftings $\rho'_2$ and $\rho'_1$.) Conversely, if $(\rho',H)$ is a normal nondegenerate representation of $\cB'$ on the Hilbert space $H$, then $E:=C_{\cB'}(\sB(G,H))$ as in \eqref{CB'def} defines a concrete von Neumann \nbd{\cB}module in $\sB(G,H)$, which gives back $\rho'$ as commutant lifting. The only critical task, nondegeneracy in Condition \ref{cvnM3}, is settled by the following result.

\bitemp[Lemma (Muhly and Solel \protect{\cite[Lemma 2.10]{MuSo02}}).~]\label{MuSolem}
If $\rho'$ is a nondegenerate normal representation of $\cB'$ on a Hilbert space $H$, then the intertwiner space $C_{\cB'}(\sB(G,H))$ acts nondegenerately on $G$.\footnote{Denote by $P\in\tMatrix{\cB&E^*\\E&\sB^a(E)}'=\Bfam{\tMatrix{b'&0\\0&\rho'(b')}\colon b'\in\cB'}$ the projection onto the invariant subspace $\cls\tMatrix{\cB&E^*\\E&\sB^a(E)}\tMatrix{G\\H}$. This subspace contains $G$ so that $P=\tMatrix{\U&0\\0&\rho'(\U)}$. Since $\rho'$ is nondegenerate, the statement follows.}
\eitemp

We find that
\bal{\label{MRfun}
(E,H)
&~\longleftrightarrow~
(\rho',H)
&
a\in\sB^a(E_1,E_2)
&~\longleftrightarrow~
a\in\sB^{bil}(H_1,H_2)
}\eal
establishes a bijective functor between the category $\cvN_\cB$ of concrete von Neumann \nbd{\cB}modules and the \hl{category $_{\cB'}\cvN$ of normal nondegenerate representations of $\cB'$} with the intertwiners $\sB^{bil}(H_1,H_2)$ as morphisms. (The preceding correspondence was established in Skeide \cite{Ske03c} as an equivalence between the category von Neumann \nbd{\cB}modules and $_{\cB'}\cvN$. As a von Neumann \nbd{\cB}module $E$, first, must be turned into a concrete von Neumann \nbd{\cB}modules $(E,E\odot G)$, the correspondence is not bijective but only an equivalence. The precise formulation above, where the functor is, really, bijective and not only an equivalence, is due to \cite{Ske06b}.)

\lf
A \hl{concrete von Neumann correspondence} over a von Neumann algebra $\cB$ is a concrete von Neumann \nbd{\cB}module $(E,H)$ with a left action of $\cB$ such that $\rho\colon\cB\rightarrow\sB^a(E)\rightarrow\sB(H)$ defines a normal (nondegenerate, of course) representation of $\cB$ on $H$. We call $\rho$ the \hl{Stinespring representation} associated with $E$.

\brem\label{Sticlrem}
The GNS-correspondence of a (normal) CP-map $T$ on $\cB$ is the unique Hilbert (von Neumann) \nbd{\cB}correspondence $E$ which has a vector $\xi\in E$ that generates $E$ as a (von Neumann) correspondence and gives back $T$ as $T(b)=\AB{\xi,b\xi}$; see Paschke \cite{Pas73}. For this GNS-correspondence $E$, the representation $\rho$ is, indeed, the Stinespring representation, while $\rho'$ is (a restriction of) the representation constructed by Arveson \cite{Arv69} in the section called ``lifting commutants''.
\erem

By $_\cB\cvN_\cB$ we denote the \hl{category of concrete von Neumann correspondences from $\cB$ to $\cB$} with the bilinear adjointable maps $a\in\sB^{a,bil}(E_1,E_2)$ as morphisms. (For adjointable maps, only left \nbd{\cB}linearity has to be checked.) We observe that $\rho(\cB)\subset\sB^a(E)=\rho'(\cB')'$, that is, $\rho'$ and $\rho$ have mutually commuting ranges. As this is very close to correspondences in the sense of Connes \cite{Con80p} (if $\cB$ is in \it{standard form}, then $\cB'\cong\cB^{op}$), we introduce the \hl{category of concrete Connes correspondences} $_\cB\ecC_\cB$ whose objects are triples $(\rho',\rho,H)$ such that $\rho'$ and $\rho$ are a pair of normal nondegenerate representations of $\cB'$ and of $\cB$, respectively, on $H$ with mutually commuting ranges, and with those maps in $\sB(H_1,H_2)$ as morphisms that intertwine both actions that of $\cB'$ and that of $\cB$. Extending the correspondence between concrete von Neumann \nbd{\cB}modules and representations of $\cB'$, we find a find bijective functor between the category of concrete von Neumann \nbd{\cB}correspondences $(E,H)$ and the category of concrete Connes correspondences $(\rho',\rho,H)$. In \cite{Ske03c} we observed this as an equivalence for von Neumann correspondences, while the bijective version for concrete von Neumann correspondences is from \cite{Ske06b}.

A last almost trivial observation (once again in \cite{Ske03c} up to equivalence and in \cite{Ske06b}, really, bijective) consists in noting that in the representation picture the roles of the representations $\rho'$ and $\rho$ are absolutely symmetric. That is, $_\cB\ecC_\cB\cong{_{\cB'}}\ecC_{\cB'}$. Therefore, if we switch $\cB$ and $\cB'$, that is, if we interprete $\rho$ as commutant lifting of $\cB$, the commutant of $\cB'$, and $\rho'$ as Stinespring representation of $\cB'$, by
\beq{\label{CBdef}
E'
~:=~
C_\cB(\sB(G,H))
~:=~
\bCB{x'\in\sB(G,H)\colon\rho(b)x'=x'b~(b\in\cB)}
}\eeq
we obtain a von Neumann \nbd{\cB'}module which is turned into a von Neumann \nbd{\cB'}correspondence by defining a left action via $\rho'$. We call $E'$ the \hl{commutant} of $E$. The commutant is a bijective functor from the category of concrete von Neumann \nbd{\cB}correspondences onto the category of concrete von Neumann \nbd{\cB'}correspondences (in each case with the bilinear adjointable maps as morphisms that are, really, the same algebra $\sB^a(E)\cap\sB^a(E')=\rho'(\cB')'\cap\rho(\cB)'$ of operators in $\sB(H)$). Obviously, $E'':=(E')'=E$.

\brem\label{gencommrem}
Muhly and Solel \cite{MuSo04} have discussed (independently) a version of the commutant for \nbd{W^*}al\-ge\-bras, called \it{\nbd{\sigma}dual}, where $\sigma$ is a faithful representation of the underlying \nbd{W^*}algebra, that must be chosen, and the \nbd{\sigma}dual depends on $\sigma$ (up to Morita equivalence of correspondences \cite{MuSo05}). An extension to correspondences from $\cA$ to $\cB$ was first done in the setting of \nbd{\sigma}duals in \cite{MuSo05}. In \cite{Ske06b} we discussed the version for von Neumann algebras and (concrete) von Neumann correspondences.

We remark that the functor $\cvN_\cB\leftrightarrow{_{\cB'}}\cvN$ in \eqref{MRfun} fits canonically into the setting of the commutant functor as $_\C\cvN_\cB\longleftrightarrow{_{\cB'}}\cvN_{\C'}$, if we consider $\C=\C'\subset\sB(\C)=\C$ as a von Neumann algebra.
\erem

The tensor product of Connes correspondences is tricky to describe in terms that do not explicitly involve the von Neumann correspondences to which they correspond. It requires that the von Neumann algebra is a \nbd{W^*}algebra in standard form and parts from Tomita-Takesaki theory and the result depends manifestly on the choice of a normal semifinite weight; see, for instance, Takesaki \cite[Section IX.3]{Tak03a}. Also the tensor product of \nbd{W^*}correspondences, although definitely less involved, still has the problem that the usual tensor product must be completed in a suitable \it{\nbd{\sigma}topology}, and this topology is defined rather \it{ad hoc}.

The tensor product two of von Neumann correspondences $E_1$ and $E_2$ is easy to obtain (and unique up to unitary equivalence): Simply construct $E_1\odot E_2\odot G$ and determine the strong closure of $E_1\,\ul{\odot}\;E_2$ in $\sB(G,E_1\odot E_2\odot G)$ or, equivalently, determine the intertwiner space $C_{\cB'}(\sB(G,E_1\odot E_2\odot G))$, a purely algebraic problem, like determining the double commutant of a \nbd{*}algebra of operators. Up to canonical isomorphism it is not important whether we construct first $E_1\,\ul{\odot}\;E_2$ and then $(E_1\,\ul{\odot}\;E_2)\odot G$ or first $E_2\odot G$ and then $E_1\odot(E_2\odot G)$. If we have concrete von Neumann correspondences $(E_1,H_1)$ and $(E_2,H_2)$ it occurs to be more adapted to construct $E_1\odot H_2$ as the space $H_2$, canonically isomorphic to $E_2\odot G$, is given from the beginning. By slight abuse of notation we shall denote the concrete von Neumann correspondence obtained in that way by $E_1\odot E_2\subset\sB(G,E_1\odot H_2)$, using the same symbol $\odot$ as for the tensor product of \nbd{C^*}correspondences. Anyway, no matter how we obtained $E_1\odot E_2\odot G$, as $(E_1\,\ul{\odot}\;E_2)\odot G$, as $E_1\odot(E_2\odot G)$ or as $E_1\odot H_2$, to fix an isomorphism from the concrete von Neumann correspondence $(E_1\odot E_2,E_1\odot E_2\odot G)$ to a concrete von Neumann correspondence $(F,K)$ simply means to fix a unitary $u\in\sB(E_1\odot E_2\odot G,K)$ that intertwines both the commutant liftings of $\cB'$ and the Stinespring representations of $\cB$.

\lf
The notations established so far allow to state and prove that the commutant establishes a bijective functor between the category $\cvN^\odot_\cB$ of product systems of concrete von Neumann \nbd{\cB}cor\-re\-spond\-ences and the category $\cvN^\odot_{\cB'}$ of product systems of concrete von Neumann \nbd{\cB'}cor\-re\-spond\-ences; see Skeide \cite{Ske03c,Ske06p3} and Muhly and Solel \cite{MuSo05p}. A morphism between two objects $E^\odot$ and $F^\odot$ in $\cvN^\odot_\cB$ is a family $a^\odot=\bfam{a_t}_{t\in\bS}$ of maps $a_t\in\sB^{a,bil}(E_t,F_t)$ that fulfills $a_s\odot a_t=a_{s+t}$ and $a_0=\id_\cB$.

We sketch this very briefly. Suppose $E^\odot$ is a product system of concrete von Neumann correspondences $(E_t,H_t)$ over $\cB$ and denote by $(\rho'_t,\rho_t,H_t)$ the corresponding concrete Connes correspondence. The familiy $u_{s,t}\colon E_s\odot E_t\rightarrow E_{s+t}$ that determines the product system structure, in the picture of Hilbert spaces is captured by unitaries $u_{s,t}\in\sB(E_s\odot H_t,H_{s+t})$ that intertwine both the actions of $\cB$ and the actions of $\cB'$ on these Hilbert spaces. The associativity condition reads $u_{r,s+t}(x_r\odot u_{s,t}(y_s\odot h_t))=u_{r+s,t}(u_{r,s}(x_s\odot y_r)\odot h_t)$. The double meaning in this formula of $u_{r,s+t}$, $u_{s,t}$, $u_{r+s,t}$ as operators between Hilbert spaces and of $u_{r,s}$ as operator between correspondences is not exactly satisfactory. We will circumvent this (purely formal) difficulty in \cite{Ske07p} by giving a different definition of the product system structure in terms of \it{representations of Hilbert modules}\footnote{This is not to be confused with the term \it{covariant representation} of a correspondence in the work of Muhly and Solel starting with \cite{MuSo98}.} where the $u_{s,t}$ will appear no longer as a defining object but as a derived one. Here we limit ourselves to explain how the product system structure of $E^\odot$ gives rise to a product system structure of $E'^\odot$ by a giving some (canonical) isomorphisms as we did in \cite{Ske03c}.

Indeed, to define a (bilinear) unitary $E'_t\odot E'_s\rightarrow E'_{t+s}$ we have to establish a unitary $u'_{t,s}\colon E'_t\odot H_s\rightarrow H_{t+s}$ intertwining the relevant representations. We have
\beqn{
E_s\odot E_t\odot G
~\cong~
E_s\odot E'_t\odot G
~\cong~
E'_t\odot E_s\odot G
~\cong~
E'_t\odot E'_s\odot G,
}\eeqn
where we used two times $E_t\odot G\cong H_t\cong E'_t\odot G$ and, in the middle, an isomorphism that, indeed, simply flips $x_s\odot y'_t\odot g$ to $y'_t\odot x_s\odot g$. (This flip is the only place where we have to compute something.) Attaching elements $b\in\cB$ and $b'\in\cB'$ in every part to the places where they act naturally, we see that the suggested isomorphism intertwines their actions. This chain of isomorphisms written down for an arbitrary tensor product of correspondences (even over different von Neumann algebras; see Remark \ref{gencommrem}) shows clearly that the commutant flips orders in tensor products:
\beq{\label{comswitch}
(E\odot F)'
~\cong~
F'\odot E';
}\eeq
see \cite[Theorem 2.3]{Ske03c} and \cite[Lemma 3.3]{MuSo05}. For our scope here, defining $u'_{t,s}$, it is sufficient to look at the chain
\beq{\label{direct}
E'_t\odot H_s
~=~
E'_t\odot\cls(E_sG)
~\cong~
E_s\odot\cls(E'_tG)
~=~
E_s\odot H_t
~\stackrel{u_{s,t}}{\cong}~
H_{s+t}.
}\eeq
This chain shows clearly how $u_{s,t}$ enters and that the flip of the elements $x'_t$ ans $y_s$ is the only thing where we are really doing something.\footnote{Behind this flip there is a sort of tensor product among a von Neumann \nbd{\cB}module and a von Neumann \nbd{\cB'}module, resulting into a von Neumann \nbd{(\cB\cap\cB')'}module. In the picture of Connes correspondences this is closely related to Sauvageot \cite{Sau80,Sau83}. We will discuss this in detail in \cite{Ske07p}.} It is routine to show that the $u'_{t,s}$ defined in that way turn $E'^\odot$ into a product system, the \hl{commutant system} of $E^\odot$.

\bex\label{commex}
Suppose $E^\otimes$ is an Arveson system. We turn $E_t$ in to the concrete von Neumann correspondence $(E_t,H_t=E_t)$ by identifying $x_t\in E_t$ with the map $\lambda\mapsto x_t\lambda$ in $\sB(\C,E_t)$. Both the Stinespring representation of $\C\subset\sB(\C)$ and the commutant lifting of $\C'=\C$ are simply the map $\lambda\mapsto\lambda\id_t$. So $E'_t=C_{\C}(\sB(\C,E_t))=E_t$. But, the unitary $u'_{t,s}\colon E'_t\otimes E'_s\rightarrow E'_{t+s}$ suggested by \eqref{direct} means: Take $x'_t\otimes y'_s$, express $y'_s\in H_s=E_s$ as $y_s1\in E_s\C$ with $y_s=y'_s$, flip the tensors to $y_s\otimes x'_t1$ and interpret $x'_t1\in E'_t\C$ as element in $x_t\in H_t=E_t=E'_t$. In this special case, thanks to $E_t=E'_t$ (a formula, that in the general case has no sense), everything can be done at once by defining $u'_{t,s}(x'_t\otimes y'_t):=u_{s,t}(y'_s\otimes x'_t)$.
\eex

As a more elaborate example we discuss now the construction from \cite{Ske03c} that, first, generalizes Arveson's construction of a product system from an \nbd{E_0}semigroup and, then, shows that the commutant of that system is the product system associated with the \nbd{E_0}semi\-group. So let $(E,H)$ be a strongly full von Neumann module over a von Neumann algebra $\cB\subset\sB(G)$. \hl{Strongly full} means that $\AB{E,E}$ generates $\cB$ as von Neumann algebra. It is easy to show (exercise!) that $E$ is strongly full, if and only if the associated commutant lifting $(\rho',H)$ is faithful. Let $\vt$ be an \nbd{E_0}semigroup of (normal unital) endomorphisms $\vt_t$ of $\sB^a(E)\subset\sB(H)$. Like Arveson we define the intertwiner space
\beqn{
E'_t
~:=~
\bCB{x'_t\in\sB(H)\colon\vt_t(a)x'_t=x'_ta~(a\in\sB^a(E))}.\footnote{Note that the interwiners are calculated in $\sB(H)$ of which $\sB^a(E)$ is a von Neumann subalgebra. Taking the analogy with \cite{Arv89} too literally would mean to determine only the interwiners in $\sB^a(E)$ what is not approximately as useful.}
}\eeqn
Since $x'^*_ty'_t\in\sB^a(E)'=\rho'(\cB')$ and $\rho'$ is faithful, we may define an inner product
\beqn{
\AB{x'_t,y'_t}
~:=~
\rho'^{-1}(x'^*_ty'_t)
}\eeqn
on $E'_t$ with values in $\cB'$. It is plain to verify that $E'_t$ with this inner product and the bimodule operation $b'_1x'_tb'_2:=\rho'(b'_1)x'_t\rho'(b'_2)$ is a von Neumann correspondence over $\cB'$. Observe that $E'_t$ acts nondegenerately on $H$ by Lemma \ref{MuSolem} and that it is the only space of intertwiners of $\vt_t$ and $\id_{\sB^a(E)}$ with this property. It follows that
\beqn{
x'_t\odot y'_s
~\mapsto~
x'_ty'_s
}\eeqn
defines an isomorphism from $E'_t\odot E'_s$ onto $E'_{t+s}$ turning $E'^\odot=\bfam{E'_t}_{t\in\bS}$ into a product system. 

Like the Arveson system of an \nbd{E_0}semigroup on $\sB(H)$, the product system $E'^\odot$ associated à la Arveson with the \nbd{E_0}semigroup $\vt$ on $\sB^a(E)\subset\sB(H)$ comes along with a \hl{faithful nondegenerate representation} $\eta'_t\colon E'_t\rightarrow\sB(H)$ on the Hilbert space $H$. By this we mean two things: Firstly, the $\eta'_t$ are isometric in the sense that $\eta'_t(x'_t)^*\eta'_t(y'_t)=\eta'_0(\AB{x'_t,y'_t})$, and the $\eta'_t$ are multiplicative in the sense that $\eta'_t(x'_t)\eta'_s(y'_s)=\eta'_{t+s}(x'_ty'_s)$. From this it follows that for each $t$ the pair $(\eta'_t,\eta'_0)$ is an isometric covariant representation of $E'_t$ in the sense of Muhly and Solel \cite{MuSo98}. In particular, $\eta'_0$ is a representation of $\cB'$. And, secondly, the $\eta'_t$ are normal and faithful in the sense that their unique extension to a representation of the \it{von Neumann linking algebra} $\sMatrix{\cB'&E'^*_t\\E'_t&\sB^a(E'_t)}$ on $\sB(H\oplus H)=M_2(\sB(H))$ is normal and faithful. One may show that this is the case, if and only if $\eta'_0$ is normal and faithful. In our case, $\eta'_t$ is simply the canonical embedding $E'_t\rightarrow\sB(H)$. In particular, $\eta'_0=\rho'$. Like for Arveson systems, speaking about a faithful nondegenerate representation is the same as speaking about a \hl{right dilation} of $E'^\odot$ to $R$ in the sense that $R$ is Hilbert space with a faithful normal representation $\rho'$ of $\cB'$ and $w'_t\colon E'_t\odot R\rightarrow R$ is a family of unitaries in $\sB^{bil}(E'_t\odot R,R)$ that iterates associatively with the product system structure. If we have such a right dilation then $\vt^{w'}_t(a)=w'_t(\id_{E'_t}\odot a)w'^*_t$ defines an \nbd{E_0}semigroup on the concrete von Neumann \nbd{\cB}module $(L,R)$ determined by $(\rho',R)$, giving back $E'^\odot$ as product system of intertwiners. Note that faithfulness of $\rho'$ implies that $L$ is strongly full and that the left action of $\cB'$ on all $E'_t$ is faithful. We say the product system $E'^\odot$ is \hl{faithful}.

Returning to the \nbd{E_0}semigroup $\vt$ on $\sB^a(E)$ and its product system \it{à la} Arveson $E'^\odot$, let us turn $E'_t$ into a concrete von Neumann correspondence $(H_t,\sigma_t,\sigma'_t)$ by defining the Hilbert space $H_t=E'_t\odot G$, the commutant lifting $\sigma_t(b)=\id_{E'_t}\odot b$ of $(\cB')'=\cB$, and the Stinespring representation $\sigma'_t(b')=b'\odot\id_G$ of $\cB'$. Then its commutant
\beqn{
E_t
~:=~
C_{\cB'}(\sB(G,H_t))
~=~
\bCB{x_t\in\sB(G,H_t)\colon\sigma'_t(b')x_t=x_tb'~(b'\in\cB')}
}\eeqn
is a concrete von Neumann correspondence over $\cB$ with left action via $\sigma_t$. If we now apply \eqref{comswitch} to $E'_t\odot H$ what we find is $H'\odot E_t=E\odot E_t$, because $(H',H)$, the commutant of $(H,H)\in{_{\cB'}}\cvN_\C$, is nothing but $(E,H)$; see Remark \ref{gencommrem}. The (bilinear!) isomorphisms $w'_t\colon E'_t\odot H\rightarrow H$ give, therefore, rise to isomorphisms $v_t:=(w'_t)'\colon E\odot E_t\rightarrow E$. (Again, we do not show that associativity is respected.) This transition between right dilations of $E'^\odot$ and left dilations of $E^\odot$ does not depend on that we started from an \nbd{E_0}semigroup on $\sB^a(E)$. In fact, we have the complete analogue of Observation \ref{lrob} including that $E^\odot$ is strongly full, if and only if $E'^\odot$ is faithful. We mention that this correspondence between left dilations of $E^\odot$ (that is, \nbd{E_0}semigroup having $E^\odot$ as associated product system) and and right dilations of $E'^\odot$ (that is, nondegenerate representations of $E'^\odot$) is due to \cite{Ske04p}. It has been generalized to not necessarily unital endomorphism semigroups and not necessarily nondegenerate representations in \cite[Theorem 3.6(3)]{Ske06p3}. Just that in \cite{Ske04p,Ske06p3} we did not yet use the terminology of dilations of a product system.

If $E^\odot$ is both strongly full and faithful, then also Observation \ref{pairob} remains true. In \cite{Ske07p} we will show that a \it{strongly continuous} product system $E^\odot$ of (concrete) von Neumann correspondences admits a \it{strongly continuous} left dilation, if it is strongly full, and a \it{strongly continuous} right dilation, if it is faithful.\footnote{Unlike continuous product systems (unital $\cB$) where existence of a unit vector in $E_1$ was automatic, it is an open problem whether strongly continuous product systems always have unit vectors. (The \nbd{C^*}proof does no longer work, because the invertibles are not open in the strong topology.) So, the results from \cite{Ske04p} (existence of a left and right dilation of the discrete subsystem) and the basic idea of \cite{Ske06} (how to turn it into a dilation of the continuous time system) become indispensable.} In this case, the unitary group according to Observation \ref{pairob} is strongly continuous and so are the two \nbd{E_0}semigroups $\vt$ and $\theta$.\footnote{If we apply this module version of Observation \ref{pairob} to $\sB^a(\cB)=\cB$, then we obtain a completely different proof of a result due to Arveson and Kishimoto \cite{ArKi92}: Every faithful normal \nbd{E_0}semigroup is the restriction of an inner automorphism group on some $\sB(H)$ to a subalgebra isomorphic to $\cB$.} This means that the commutant of $E^\odot$ is also derived from a strongly continuous \nbd{E_0}semigroup and, therefore, possesses a \hl{strongly continuous} structure. We do not describe here the definition of \it{strongly continuous} product systems. Apart from missing space, at the time being we have more than one candidate for a definition, and all candidates work well. We did not yet find out which one we should consider the best one. Anyway, the results will allow to show the following theorem.

\bthm
The commutant of a strongly continuous, strongly full and faithful product system is strongly continuous, strongly full and faithful, too.
\ethm

Muhly and Solel \cite{MuSo05p} have a similar result for \it{measurable} product systems under separability assumptions. However, while our proof relies essentially on the product system structure (in that we have to construct a left and a right dilation and to use the semigroup structure encoded by them), their proof is rather a result on general measurable bundles of correspondences and works by a reduction to the analogue result for von Neumann algebras due to Effros \cite{Eff65}. In general it is far from being obvious why a bundle of intertwiners between bundles of Banach modules should admit (strongly) continuous sections.

\brem\label{finrem}
The list of dualities may be extended. For instance, the fact that (by using \it{quasi orthonormal bases} of von Neumann \nbd{\cB}modules, as suitable substitute for orthonormal bases of Hilbert spaces) every von Neumann \nbd{\cB}module is a complemented submodule of a \it{free} von Neumann \nbd{\cB}module, may be used to prove the \it{amplification-induction theorem} on the representations $\rho'$ of $\cB'$. In the presence of invariant vector states there is a duality between CP-maps from $\cA$ to $\cB$ and CP-maps from $\cB'$ to $\cA'$ (Albeverio and Hoegh-Krohn \cite{AlHK78}) that includes a duality between \it{tensor dilations} of a CP-maps on $\cB$ and \it{extensions} from $\cB'$ to $\sB(G)$ of the dual of that CP-map; see Gohm and Skeide \cite{GoSk05}. Applying the duality of CP-maps to the canonical embedding of a subalgebra $\cA\subset\cB$ into $\cB$ (both in standard form) and translating back the dual map $\cB'\rightarrow\cA'$ into a map $\cB\rightarrow\cA$ via twofold \it{Tomita conjugation}, one obtains the \it{Accardi-Cecchini conditional expectation} \cite{AcCe82} that coincides with the usual conditional expectation whenever the latter exists; see also Accardi and Longo \cite{Lon84,AcLo93}.

We mention the duality between Rieffel's \it{Eilenberg-Watts theorem} \cite{Rie74a} about functors between categories of representations of von Neumann algebras and Blecher's \it{Eilenberg-Watts theorem} \cite{Ble97} about functors between categories of Hilbert modules. When the latter is restricted to von Neumann modules, the two \it{Eilenberg-Watts theorems} are dual to each other under the commutant; see \cite{Ske06b}.

Last but surely not least there is the duality between the product system of a CP-semigroup in Bhat and Skeide \cite{BhSk00} and the product system constructed from the same CP-semigroup by Muhly and Solel \cite{MuSo02}. The latter is the the commutant of the former, a problem left open in \cite{MuSo02} that lead to the notion of commutant of correspondences and product systems in \cite{Ske03c}.
\erem

%\newpage

\setlength{\baselineskip}{2.5ex}

%\renewcommand{\tt}[1]{\texttt{\small #1}}
%\bibliography{mybib}

\newcommand{\Swap}[2]{#2#1}\newcommand{\Sort}[1]{}
\providecommand{\bysame}{\leavevmode\hbox to3em{\hrulefill}\thinspace}
\providecommand{\MR}{\relax\ifhmode\unskip\space\fi MR }
\providecommand{\MRhref}[2]{%
  \href{http://www.ams.org/mathscinet-getitem?mr=#1}{#2}
}
\providecommand{\href}[2]{#2}

%\listofOWs

\end{document}